# THE STABILITY OF CONDITIONAL MARKOV PROCESSES AND MARKOV CHAINS IN RANDOM ENVIRONMENTS

By Ramon van Handel

*Princeton University*

We consider a discrete time hidden Markov model where the signal is a stationary Markov chain. When conditioned on the observations, the signal is a Markov chain in a random environment under the conditional measure. It is shown that this conditional signal is weakly ergodic when the signal is ergodic and the observations are nondegenerate. This permits a delicate exchange of the intersection and supremum of $\sigma$-fields, which is key for the stability of the non-linear filter and partially resolves a long-standing gap in the proof of a result of Kunita [*J. Multivariate Anal.* **1** (1971) 365–393]. A similar result is obtained also in the continuous time setting. The proofs are based on an ergodic theorem for Markov chains in random environments in a general state space.

**1. Introduction.** Consider a discrete time Markov chain $(X_n)_{n\in\mathbb{Z}_+}$ and a random process $(Y_n)_{n\in\mathbb{Z}_+}$ such that $Y_n$ and $Y_m$ $(n\neq m)$ are conditionally independent given $(X_n)_{n\in\mathbb{Z}_+}$ and such that the conditional distribution of $Y_n$ given $(X_n)_{n\in\mathbb{Z}_+}$ depends only on $X_n$. Then the pair $(X_n, Y_n)_{n\in\mathbb{Z}_+}$ defines a hidden Markov model, where the observation process $(Y_n)_{n\in\mathbb{Z}_+}$ provides indirect information on the signal process $(X_n)_{n\in\mathbb{Z}_+}$. Models of this form have a wide array of applications in statistics, engineering and finance, and possess a rich theory of statistical inference [7]. Of particular interest in the present paper is the filtering problem, which aims to estimate the current state $X_n$ of the signal given the observation history $(Y_k)_{0\leq k\leq n}$ by computing the regular conditional probability $\mathbf{P}(X_n\in\cdot\,|(Y_k)_{0\leq k\leq n})$. A similar class of problems can also be formulated in continuous time.

This paper is concerned with the long time properties of the nonlinear filter, that is, we are interested in the behavior of the regular conditional

---

probabilities $\Pi_n = \mathbf{P}(X_n \in \cdot | (Y_k)_{0 \le k \le n})$ as $n \to \infty$, in the case that the signal possesses an invariant probability measure $\pi$. The investigation of such problems in general hidden Markov models has a long history, starting with the pioneering work of Kunita [23] (in the continuous time setting) on the stationary behavior of the mean square estimation error of the nonlinear filter. To study this problem, he established the following key result [23], Theorem 3.3: for any invariant measure $\pi$ of the signal, the filtering process $(\Pi_n)_{n \in \mathbb{Z}_+}$ possesses a unique invariant measure with barycenter $\pi$ if and only if the signal is ergodic in a particular sense (see below).

A different but closely related problem of interest is the stability of nonlinear filters. Denote by $\mathbf{P}^\mu$ the law of $(X_n, Y_n)_{n \in \mathbb{Z}_+}$ with the initial law $X_0 \sim \mu$, and write the corresponding filter as $\Pi_n^\mu = \mathbf{P}^\mu(X_n \in \cdot | (Y_k)_{0 \le k \le n})$. In practice, the initial measure $\mu$ (the Bayesian prior) is rarely known precisely, and it is thus highly desirable that the filter $\Pi_n^\mu$ becomes insensitive to the choice of $\mu$ as $n \to \infty$ (e.g., as in Theorem 5.2 below). When this is the case, the filter is said to be stable. In a pioneering paper, Ocone and Pardoux [25] used Kunita's theorem to establish that stability of the filter is inherited from the ergodicity of the signal process.

The asymptotic properties of nonlinear filters have received considerable attention in recent years (see, e.g., [12] and the references therein). Beside the fundamental interest of the topic, results in this direction have a variety of applications, which include uniform convergence of filter approximations [5, 6, 13, 14], maximum likelihood estimation [7, 8, 19], stochastic control [18, 31] and estimation error bounds [3, 23]. In various specific cases one can even obtain detailed quantitative information about the rate of stability of the filter (see [12] for references). In the general setting, however, little is known about the asymptotic properties of nonlinear filters beyond the work of Kunita [23] and subsequent papers, such as [25], which rely directly on the approach of [23] (but see [36]).

Unfortunately, as was pointed out in [1], there is a serious gap in the proof of the main result in [23]. To describe the problem, let us suppose that the signal process possesses an invariant probability measure $\pi$. Then $\mathbf{P}^\pi$ is a stationary measure, and we can therefore extend the stationary hidden Markov model to two-sided time $(X_n, Y_n)_{n \in \mathbb{Z}}$ by a standard argument. Denote by $\mathbf{P}$ the extension of $\mathbf{P}^\pi$ to two-sided time, and define the $\sigma$-fields $\mathcal{F}_I^X = \sigma\{X_n : n \in I\}$ and $\mathcal{F}_I^Y = \sigma\{Y_n : n \in I\}$ ($I \subset \mathbb{Z}$). The key step in Kunita's proof is to argue that his result would follow if we could establish that the following identity holds true:

$$\bigcap_{n \ge 0} \mathcal{F}_{]-\infty, 0]}^Y \vee \mathcal{F}_{]-\infty, -n]}^X = \mathcal{F}_{]-\infty, 0]}^Y \qquad \mathbf{P}\text{-a.s.}$$

He proceeds to argue as follows. Suppose that the signal satisfies the following ergodicity condition: $\bigcap_{n \geq 0} \mathcal{F}^X_{]-\infty, -n]}$ is $\mathbf{P}$-a.s. trivial. Then

$$\bigcap_{n \geq 0} \mathcal{F}^Y_{]-\infty, 0]} \vee \mathcal{F}^X_{]-\infty, -n]} \overset{?}{=} \mathcal{F}^Y_{]-\infty, 0]} \vee \bigcap_{n \geq 0} \mathcal{F}^X_{]-\infty, -n]} = \mathcal{F}^Y_{]-\infty, 0]} \qquad \mathbf{P}\text{-a.s.}$$

The exchange of the intersection and supremum of $\sigma$-fields is not at all obvious, however, and no proof of this assertion is provided in [23]. Indeed, this exchange is not permitted in general, as an illuminating counterexample in [1] shows.

It is important to note, on the other hand, that all known counterexamples rely in an essential way on the degeneracy of the observation model, that is, $Y_k = h(X_k)$ for some function $h$ without any additional noise. It is therefore tempting to conjecture that the exchange of intersection and supremum is always permitted provided that the observations are nondegenerate, which is most naturally imposed in our general setting by requiring that the conditional law of $Y_n$ given $(X_k)_{k \in \mathbb{Z}}$ satisfies

$$\mathbf{P}(Y_n \in A | (X_k)_{k \in \mathbb{Z}}) = \int I_A(du) g(X_n, u) \varphi(du) \qquad \mathbf{P}\text{-a.s.,}$$

where $\varphi$ is a fixed reference measure and $g$ is a strictly positive function. Though no counterexamples are known, it is unclear whether or not this is the case, and the (positive or negative) verification of this conjecture remains an open problem.

From the work of Budhiraja [4] and of Baxendale, Chigansky and Liptser [1], and from the results of Section 5 below, it is clear that Kunita's exchange of intersection and supremum and its time-reversed cousin

$$\bigcap_{n \geq 0} \mathcal{F}^Y_{]-\infty, 0]} \vee \mathcal{F}^X_{]-\infty, -n]} \overset{?}{=} \mathcal{F}^Y_{]-\infty, 0]}$$

and

$$\bigcap_{n \geq 0} \mathcal{F}^Y_{[0, \infty[} \vee \mathcal{F}^X_{[n, \infty[} \overset{?}{=} \mathcal{F}^Y_{[0, \infty[} \qquad \mathbf{P}\text{-a.s.}$$

lie at the heart of the qualitative asymptotic theory of nonlinear filtering. The main result of this paper, Theorem 4.2, establishes that both these identities do indeed hold under conditions that are only mildly stronger than those assumed by Kunita. Given an invariant probability measure $\pi$ of the signal process, we assume the following:

1. The signal is ergodic in the following sense:

$$\|\mathbf{P}^{\delta_z}(X_n \in \cdot) - \pi\|_{\mathrm{TV}} \xrightarrow{n \to \infty} 0 \qquad \text{for } \pi\text{-a.e. } z,$$

where $\| \cdot \|_{\mathrm{TV}}$ is the total variation norm (Assumption 3.1 below).

2. The observations are nondegenerate (Assumption 3.2 below).

These assumptions are satisfied by the vast majority of stationary hidden Markov models of practical interest, including the important case of aperiodic and positive Harris recurrent signals with nondegenerate observations. Note that we do not require the Feller assumption, and that we allow for signal and observation processes with arbitrary Polish state spaces (the Polish assumption guarantees an abundance of regular conditional probabilities). The latter has the additional advantage that our results extend directly to the continuous time setting (Section 6).

Beside our main result, this paper contains two additional results which are of independent interest. First, as we will discuss shortly, the proof of our main result hinges on the ergodic theory of Markov chains in random environments as developed by Cogburn [10, 11] and Orey [26] for countable state spaces. In Section 2, we prove the counterpart of a result from [11] for Markov chains in random environments on general Polish state spaces (Theorem 2.3). This result is not specific to hidden Markov models, and could be relevant in other settings.

Second, we will show in Section 5 that the permissibility of the exchange of intersection and supremum leads to the stability of the nonlinear filter in a much stronger sense than was previously established in [1, 4, 25]. A special case of our main stability theorem (Theorem 5.2) is the following result: if the signal is aperiodic and positive Harris recurrent, and if the observations are nondegenerate, then

$$\|\Pi_n^\mu - \Pi_n^\nu\|_{\mathrm{TV}} \xrightarrow{n \to \infty} 0 \qquad \mathbf{P}^\gamma\text{-a.s. for all } \mu, \nu, \gamma.$$

Similar results hold in the continuous time setting (Section 6).

The remainder of this section is devoted to a guided tour through our proofs.

1.1. *The method of von Weizsäcker and the conditional signal.* In [37], von Weizsäcker has studied the exchange of intersection and supremum problem in a general setting. Following his approach, one can establish the following illuminating result. Let $\mathcal{G}_n$, $n \in \mathbb{N}$ be a decreasing family of countably generated $\sigma$-fields and let $\mathcal{F}$ be another countably generated $\sigma$-field. Then

$$\bigcap_{n \in \mathbb{N}} \mathcal{F} \vee \mathcal{G}_n = \mathcal{F} \quad \mathbf{P}\text{-a.s.} \quad \text{iff} \quad \bigcap_{n \in \mathbb{N}} \mathcal{G}_n \text{ is } P^{\mathcal{F}}(\omega, \cdot)\text{-a.s. trivial} \quad \text{for } \mathbf{P}\text{-a.e. } \omega,$$

where $P^{\mathcal{F}}(\omega, \cdot)$ is a version of the regular conditional probability $\mathbf{P}(\cdot|\mathcal{F})$. It would appear at first glance that $\mathbf{P}$-a.s. triviality of the tail $\sigma$-field $\bigcap_{n \in \mathbb{N}} \mathcal{G}_n$ automatically implies that it is also $\mathbf{P}(\cdot|\mathcal{F})$-a.s. trivial; after all, it is elementary that $\mathbf{P}(A|\mathcal{F}) = \mathbf{P}(A)$ $\mathbf{P}$-a.s. whenever $\mathbf{P}(A) = 0$ or $\mathbf{P}(A) = 1$. However,

the tail $\sigma$-field is not countably generated, so we cannot eliminate the dependence of the exceptional set on $A$. Verification of $\mathbf{P}(\cdot|\mathcal{F})$-a.s. triviality is thus a nontrivial problem.

Despite its generality, the result of von Weizsäcker is rarely used in the literature. In many cases the result is difficult to apply, as a tractable characterization of the conditional measure $\mathbf{P}(\cdot|\mathcal{F})$ is typically not available. In our setting, however, a fortuitous observation makes this approach much more attractive: when conditioned on the observations, the signal process remains an (albeit nonhomogeneous) Markov process whose transition probabilities depend on the observed sample path of the observation process. This observation dates back to the work of Stratonovich [33], and has recently been applied to obtain quantitative stability results for various special filtering models [7, 20, 35]. In these references a time horizon $N$ is fixed and the signal is considered under the conditional measure $\mathbf{P}(\cdot|\mathcal{F}^Y_{[0,N]})$, while we will work under the conditional measure $\mathbf{P}(\cdot|\mathcal{F}^Y_{[0,\infty[})$, but this difference does not affect the Markov property of the conditional signal.

Our basic strategy is thus as follows. Note that by the above discussion

$$\bigcap_{n\geq 0} \mathcal{F}^Y_{[0,\infty[} \vee \mathcal{F}^X_{[n,\infty[} = \mathcal{F}^Y_{[0,\infty[} \qquad \mathbf{P}\text{-a.s.}$$

would be established if we could show that

$$\mathcal{T}^X = \bigcap_{n\geq 0} \mathcal{F}^X_{[n,\infty[} \text{ is } \mathbf{P}(\cdot|\mathcal{F}^Y_{[0,\infty[})\text{-a.s. trivial} \qquad \mathbf{P}\text{-a.s.}$$

We therefore aim to show that the signal $(X_n)_{n\geq 0}$, which is a nonhomogeneous Markov process under the regular conditional probability $\mathbf{P}(\cdot|\mathcal{F}^Y_{[0,\infty[})$, has trivial tail $\sigma$-field $\mathcal{T}^X$ for almost every observation path, provided our ergodicity and nondegeneracy assumptions are satisfied. The time-reversed result follows similarly.

1.2. *Markov chains in random environments.* To obtain our main result, we must now show that tail triviality of the signal process under the conditional measure is inherited from the ergodicity of the signal process under the original probability measure. In the following, we will often refer to the signal process under the conditional measure as the conditional signal.

To fix some ideas, consider the case of a time homogeneous finite state Markov chain. In this setting, ergodicity (and hence tail triviality) is determined entirely by the graph of the chain, and not by the precise values of the transition probabilities. In particular, for one such chain to inherit ergodicity from another chain, it suffices that the two chains have the same graph, or, in probabilistic terms, that their transition probabilities are mutually absolutely continuous. That a similar statement holds in a general state space can be inferred, for example, from [28], Theorem 2.1.

The problem in our setting is that the conditional signal is not time homogeneous. Nonetheless, the transition probability of the conditional signal $K_n(x, \cdot) = \mathbf{P}(X_n \in \cdot | X_{n-1} = x, \mathcal{F}_{[0,\infty[}^Y)$ satisfies a key homogeneity property: it is easily seen [using the stationarity of $\mathbf{P}$ and the Markov property of $(X_n, Y_n)_{n \in \mathbb{Z}}$] that $n \mapsto K_n$ is a stationary stochastic process. The conditional signal is thus a Markov chain in a random environment in the terminology of Cogburn, who established [11], Section 3, that the ergodicity of such a process in a finite (or countable) state space is determined by its graph in essentially the same manner as for time homogeneous chains. This suggests that to prove our result, it suffices to show that the transition probabilities of the conditional signal and of the signal are equivalent.

As is perhaps to be expected, things are not quite so straightforward in practice. First, even in a finite state space, the conditional signal under $\mathbf{P}(\cdot | \mathcal{F}_{[0,\infty[}^Y)$ does not fit in the framework of Cogburn as the ergodic theory of Markov chains in random environments relies on the availability of all environmental variables $(Y_k)_{k \in \mathbb{Z}}$. In order to apply the result of Cogburn, we must therefore condition not on $\mathcal{F}_{[0,\infty[}^Y$ but on $\mathcal{F}_{\mathbb{Z}}^Y$. It is then necessary to establish two things: that

$$\mathbf{P}(X_n \in \cdot | X_{n-1} = i, \mathcal{F}_{\mathbb{Z}}^Y) \sim \mathbf{P}(X_n \in \cdot | X_{n-1} = i) \qquad \text{for all } i \ \mathbf{P}\text{-a.s.,}$$

so that the ergodicity of the signal process under $\mathbf{P}$ implies the ergodicity of the signal process under $\mathbf{P}(\cdot | \mathcal{F}_{\mathbb{Z}}^Y)$ by the result of Cogburn, and that

$$\mathbf{P}((X_n)_{n \geq 0} \in \cdot | \mathcal{F}_{\mathbb{Z}}^Y) \sim \mathbf{P}((X_n)_{n \geq 0} \in \cdot | \mathcal{F}_{[0,\infty[}^Y) \qquad \mathbf{P}\text{-a.s.,}$$

so that triviality of $\mathcal{T}^X$ under $\mathbf{P}(\cdot | \mathcal{F}_{\mathbb{Z}}^Y)$ implies triviality of $\mathcal{T}^X$ under $\mathbf{P}(\cdot | \mathcal{F}_{[0,\infty[}^Y)$. We will prove these identities in Sections 3 and 4 using a coupling argument; it is here that the nondegeneracy of the observations is required. Once these facts have been established, von Weizsäcker's argument completes the proof.

Unlike Cogburn's results, however, our results are not restricted to finite or countable state spaces. Our first order of business is therefore to extend the necessary result from [11] to the setting of general Polish state spaces. As with ordinary Markov chains in general state spaces, the general case requires significantly more sophisticated tools than are needed in the countable setting. Our general result in Section 2 is inspired by the elegant martingale methods of Derriennic [16] and of Papangelou [28] for ordinary Markov chains in general state spaces.

1.3. *Organization of the paper.* This paper is organized as follows.

In Section 2, we introduce the general model for a Markov chain in a random environment. The main result, Theorem 2.3, establishes that weak ergodicity, tail triviality and irreducibility are equivalent for stationary Markov chains in random environments. This result is key for the proof of our main result.

In Section 3, we introduce the general hidden Markov model. We begin by proving that this model fits in the framework of Section 2 if we condition on the complete observation record $(Y_n)_{n \in \mathbb{Z}}$ (Lemma 3.3). The main result of this section, Theorem 3.4, establishes that the conditional signal is ergodic provided that the ergodicity and nondegeneracy Assumptions 3.1 and 3.2 are satisfied. The proof proceeds in two steps. First, we show that the result would follow from ergodicity of the signal and the equivalence of the conditional and unconditional transition probabilities (Lemma 3.5). Next, we show that this equivalence does in fact hold if we additionally assume nondegenerate observations (Lemma 3.8). Of independent interest is Lemma 3.7, which is used repeatedly in the following sections.

In Section 4, we complete the proof of the main result of this paper (Theorem 4.2). First, we develop the argument of von Weizsäcker in our setting (Section 4.1). The remainder of the section is devoted to proving that $\mathbf{P}((X_n)_{n \geq 0} \in \cdot | \mathcal{F}_{\mathbb{Z}}^Y) \sim \mathbf{P}((X_n)_{n \geq 0} \in \cdot | \mathcal{F}_{[0,\infty[}^Y)$ $\mathbf{P}$-a.s. (the relevance of which was discussed above).

Section 5 establishes that our main result implies stability of the filter (Theorem 5.2). The key connection between Theorems 5.2 and 4.2 is the expression in Lemma 5.6 for the Radon–Nikodym derivative between differently initialized filters.

In Section 6, we extend our main results to the continuous time setting.

Finally, Section 7 contains a brief discussion on the implications of our main result for the gap in the result of Kunita [23].

## 2. Markov chains in random environments.

2.1. *The canonical setup and main result.* Throughout this paper, we operate in the following canonical setup. We consider the pair $(X_n, Y_n)_{n \in \mathbb{Z}}$, where $X_n$ takes values in the Polish space $E$ and $Y_n$ takes values in the Polish space $F$. We realize these processes on the canonical path space $\Omega = \Omega^X \times \Omega^Y$ with $\Omega^X = E^{\mathbb{Z}}$ and $\Omega^Y = F^{\mathbb{Z}}$, such that $X_n(x, y) = x(n)$ and $Y_n(x, y) = y(n)$. Denote by $\mathcal{F}$ the Borel $\sigma$-field on $\Omega$, and introduce the natural filtrations

$$\mathcal{F}_n^X = \sigma\{X_k : k \leq n\}, \qquad \mathcal{F}_n^Y = \sigma\{Y_k : k \leq n\}, \qquad \mathcal{F}_n = \mathcal{F}_n^X \vee \mathcal{F}_n^Y$$

for $n \in \mathbb{Z}$, as well as the $\sigma$-fields

$$\mathcal{F}_I^X = \sigma\{X_k : k \in I\}, \qquad \mathcal{F}_I^Y = \sigma\{Y_k : k \in I\}, \qquad \mathcal{F}_I = \mathcal{F}_I^X \vee \mathcal{F}_I^Y$$

for $I \subset \mathbb{Z}$. For simplicity of notation, we set

$$\mathcal{F}^X = \mathcal{F}_{\mathbb{Z}}^X, \qquad \mathcal{F}^Y = \mathcal{F}_{\mathbb{Z}}^Y, \qquad \mathcal{F}_+^X = \mathcal{F}_{[0,\infty[}^X, \qquad \mathcal{F}_+^Y = \mathcal{F}_{[0,\infty[}^Y$$

and we will denote by $Y$ the $F^{\mathbb{Z}}$-valued random variable $(Y_k)_{k \in \mathbb{Z}}$. The canonical shift $\Theta \colon \Omega \to \Omega$ is defined as $\Theta(x, y)(m) = (x(m+1), y(m+1))$.

In the following sections we will introduce a measure on $(\Omega, \mathcal{F})$ which defines a hidden Markov model. In the present section, however, it will be more convenient to attach a somewhat different interpretation to our canonical setup. To this end, consider a probability kernel of the form $P^X \colon E \times \Omega^Y \times \mathcal{B}(E) \to [0, 1]$, where $\mathcal{B}(E)$ denotes the Borel $\sigma$-field of $E$. We will define a stationary probability measure $\mathbf{P}$ on $(\Omega, \mathcal{F})$ such that the following holds a.s. for every $n \in \mathbb{Z}$:

$$\mathbf{P}(X_{n+1} \in A | \mathcal{F}_n^X \vee \mathcal{F}^Y) = P^X(X_n, Y \circ \Theta^n, A).$$

Then $X_n$ is interpreted as a Markov chain in a random environment: the environment is the sequence $Y$, and $X_n$ is a nonhomogeneous Markov process, for almost every path $Y$, under the regular conditional probability $\mathbf{P}(\cdot | \mathcal{F}^Y)$.

REMARK 2.1. Markov chains in random environments were studied extensively by Cogburn [10, 11] and by Orey [26] in the case that $E$ is countable. The purpose of this section is to extend a result in [11] to the general setting in which $E$ is Polish. It should be noted that in these papers the kernel $P^X(x, y, A)$ is assumed to depend only on $y(0)$, rather than on the entire path $y = (y(k))_{k \in \mathbb{Z}}$. This difference is immaterial, however, and the current notation fits particularly well with the hidden Markov model which will be studied in the rest of the paper.

We proceed to construct $\mathbf{P}$. Our model consists of three ingredients:

1. The probability kernel $P^X \colon E \times \Omega^Y \times \mathcal{B}(E) \to [0, 1]$.
2. A probability kernel $\mu \colon \Omega^Y \times \mathcal{B}(E) \to [0, 1]$ such that

$$\int P^X(z, y, A) \mu(y, dz) = \mu(\Theta y, A) \qquad \text{for all } y \in \Omega^Y, A \in \mathcal{B}(E).$$

3. A probability measure $\mathbf{P}^Y$ on $(\Omega^Y, \mathcal{F}^Y)$ which is invariant under the shift, that is, $\mathbf{P}^Y(Y \in A) = \mathbf{P}^Y(Y \circ \Theta \in A)$ for all $A \in \mathcal{F}^Y$.

For every $n \in \mathbb{N}$, define the probability kernel $\mathbf{P}^{(n)} \colon \Omega^Y \times \mathcal{F}_{[-n,n]}^X \to [0, 1]$ as

$$\mathbf{P}_y^{(n)}(A) = \int I_A(x) P^X(x(n-1), \Theta^{n-1} y, dx(n)) \cdots$$
$$\times P^X(x(-n), \Theta^{-n} y, dx(-n+1)) \mu(\Theta^{-n} y, dx(-n)).$$

Then $\mathbf{P}_y^{(n+1)}|_{\mathcal{F}_{[-n,n]}^X} = \mathbf{P}_y^{(n)}$, so that we can define a probability kernel

$$\mathbf{P}_\cdot \colon \Omega^Y \times \mathcal{F}^X \to [0, 1], \qquad \mathbf{P}_y|_{\mathcal{F}_{[-n,n]}^X} = \mathbf{P}_y^{(n)} \qquad \text{for all } n, y$$

by the usual Kolmogorov extension argument. We now define the probability measure $\mathbf{P}$ on $(\Omega, \mathcal{F})$ by setting

$$\mathbf{P}(A) = \int I_A(x, y) \mathbf{P}_y(dx) \mathbf{P}^Y(dy) \qquad \text{for all } A \in \mathcal{F}.$$

In addition to the probability measure $\mathbf{P}$ and the kernel $\mathbf{P}_y$, we introduce a probability kernel $\mathbf{P}_{\cdot,\cdot} : E \times \Omega^Y \times \mathcal{F}_+^X \to [0, 1]$ by setting for $A \in \mathcal{F}_{[0,n]}^X$

$$\mathbf{P}_{z,y}(A) = \int I_A(x) P^X(x(n-1), \Theta^{n-1} y, dx(n)) \cdots$$
$$\times P^X(x(1), \Theta y, dx(2)) P^X(x(0), y, dx(1)) \delta_z(dx(0)),$$

where $\delta_z(A) = I_A(z)$, and again extending by the Kolmogorov extension argument. The following is an easy consequence of our definitions.

LEMMA 2.2. *The following properties hold true:*

1. *The following holds for all $A \in \mathcal{F}_+^X$, $z \in E$, $y \in \Omega^Y$:*

$$\mathbf{E}_{z,y}(I_A \circ \Theta) = \int P^X(z, y, dz') \mathbf{P}_{z', \Theta y}(A).$$

2. $\mathbf{P}_{\Theta y}(A) = \mathbf{E}_y(I_A \circ \Theta)$ *for all $y \in \Omega^Y$, $A \in \mathcal{F}^X$.*
3. $\mathbf{P}$ *is invariant under the shift $\Theta : \Omega \to \Omega$, that is, $\mathbf{P}((X_k, Y_k)_{k \in \mathbb{Z}} \in A) = \mathbf{P}((X_{k+n}, Y_{k+n})_{k \in \mathbb{Z}} \in A)$ for all $A \in \mathcal{F}$, $n \in \mathbb{Z}$.*
4. *The following hold $\mathbf{P}$-a.s. for $A \in \mathcal{F}^X$, $B \in \mathcal{F}_+^X$, $n \in \mathbb{Z}$:*

$$\mathbf{E}(I_A \circ \Theta^n | \mathcal{F}^Y) = \mathbf{P}_{Y \circ \Theta^n}(A), \mathbf{E}(I_B \circ \Theta^n | \mathcal{F}_n^X \vee \mathcal{F}^Y) = \mathbf{P}_{X_n, Y \circ \Theta^n}(B).$$

PROOF. Elementary. $\square$

The goal of this section is to prove the following theorem. In the case that $E$ is countable, a similar result can be found in [11], Section 3.

THEOREM 2.3. *The following are equivalent.*

1. $\|\mathbf{P}_{z,y}(X_n \in \cdot) - \mathbf{P}_{z',y}(X_n \in \cdot)\|_{\text{TV}} \xrightarrow{n \to \infty} 0$ *for $(\mu \otimes \mu) \mathbf{P}^Y$-a.e. $(z, z', y)$.*
2. *The tail $\sigma$-field $\mathcal{T}^X = \bigcap_{n \geq 0} \mathcal{F}_{[n,\infty[}^X$ is a.s. trivial in the following sense:*

$$\mathbf{P}_{z,y}(A) = \mathbf{P}_{z,y}(A)^2 = \mathbf{P}_{z',y}(A) \qquad \text{for all } A \in \mathcal{T}^X \text{ and } (z, z', y) \in H,$$

*where $H$ is a fixed set (independent of $A$) of $(\mu \otimes \mu) \mathbf{P}^Y$-full measure.*
3. *For $(\mu \otimes \mu) \mathbf{P}^Y$-a.e. $(z, z', y)$, there is an $n \in \mathbb{N}$ such that the measures $\mathbf{P}_{z,y}(X_n \in \cdot)$ and $\mathbf{P}_{z',y}(X_n \in \cdot)$ are not mutually singular.*

When the first condition of this theorem holds, the Markov chain in the random environment is said to be *weakly ergodic*; when the second condition holds, it is said to be *tail trivial*; and when the last condition holds, it is said to be *irreducible*. Our goal is to prove that these notions are equivalent.

2.2. *Proof of Theorem 2.3.* The implication $1 \Rightarrow 3$ of Theorem 2.3 is trivial; thus, it suffices to show that $2 \Rightarrow 1$ and $3 \Rightarrow 1, 2$. Our approach below is partially inspired by the martingale methods of Derriennic [16] and of Papangelou [28] for ordinary Markov chains in general state spaces, and by the work of Cogburn [11] for countable Markov chains in random environments.

We begin by stating two preliminary lemmas which are in essence well-known results. The first lemma below shows that the total variation norm of a kernel is a measurable function; the second lemma shows that $2 \Rightarrow 1$ in Theorem 2.3.

LEMMA 2.4. *Let $(G, \mathcal{G})$ be a measurable space, $(K, \mathcal{K})$ be a measurable space with $\mathcal{K}$ a countably generated $\sigma$-field, and $\rho: G \times \mathcal{K} \to \mathbb{R}$ be a finite kernel. Then the map $g \mapsto \|\rho(g, \cdot)\|_{\mathrm{TV}}$ is measurable.*

PROOF. As $\mathcal{K}$ is countably generated, there is a sequence $\{I_n\}$ of refining partitions $I_n = \{E_1^n, \ldots, E_n^n\}$ of $K$ such that $\mathcal{K} = \sigma\{I_n : n \in \mathbb{N}\}$. But then

$$\sum_{k=1}^{n} |\rho(g, E_k^n)| = \|\rho(g, \cdot)|_{\sigma\{I_n\}}\|_{\mathrm{TV}} \nearrow \|\rho(g, \cdot)\|_{\mathrm{TV}} \qquad \text{as } n \to \infty$$

for all $g \in G$ (see, e.g., [27], page 1635). As $g \mapsto \rho(g, E_k^n)$ is measurable for every $k, n$, the above limit is also measurable and the result follows. □

The proof of the following result follows closely along the lines of the proof of [29], Proposition 6.2.4, and is therefore omitted.

LEMMA 2.5. *Let $H$ be a set of $(\mu \otimes \mu) \mathbf{P}^Y$-full measure. If*

$$\mathbf{P}_{z,y}(A) = \mathbf{P}_{z,y}(A)^2 = \mathbf{P}_{z',y}(A) \qquad \text{for all } A \in \mathcal{T}^X \text{ and } (z, z', y) \in H,$$

*then $\|\mathbf{P}_{z,y}(X_n \in \cdot) - \mathbf{P}_{z',y}(X_n \in \cdot)\|_{\mathrm{TV}} \xrightarrow{n \to \infty} 0$ for all $(z, z', y) \in H$. In particular, if condition 2 of Theorem 2.3 holds, then so does condition 1.*

Before we proceed, we state an additional lemma on general Markov chains which will be used several times. The construction of the set $H$ below follows closely along the lines of [27], pages 1636–1637, so the proof is omitted.

LEMMA 2.6. *Let $\mathbf{P}^z$ be the law of a Markov process $(Z_k)_{k \geq 0}$ given $Z_0 = z$, and let $\nu$ be a stationary probability for this Markov process. Then for any set $\tilde{H}$ of $\nu$-full measure, there is a subset $H \subset \tilde{H}$ of $\nu$-full measure such that*

$$\mathbf{P}^z(Z_n \in H \text{ for all } n \geq 0) = 1 \qquad \text{for all } z \in H.$$

We now proceed with the proof of Theorem 2.3. Let us introduce certain *skew* Markov chains which will be useful in what follows. Define $U_n = (X_n, Y \circ \Theta^n)$; then evidently $U_n$ is an $E \times \Omega^Y$-valued stationary Markov chain under $\mathbf{P}$, whose stationary measure $\lambda(A) = \mathbf{P}(U_n \in A)$ for all $n \in \mathbb{Z}$, $A \in \mathcal{B}(E \times \Omega^Y)$ and transition probability kernel $P^U : E \times \Omega^Y \times \mathcal{B}(E \times \Omega^Y) \to [0, 1]$ are given by

$$\lambda(A) = \int I_A(z, y) \mu(y, dz) \mathbf{P}^Y(dy), \qquad P^U(z, y, B \times C) = P^X(z, y, B) I_C(\Theta y),$$

while $U_n$ is a Markov process with the same transition probability kernel $P^U$ but with the initial measures $\delta_{z,y}$ and $\mu(y, \cdot)$ under $\mathbf{P}_{z,y}$ and $\mathbf{P}_y$, respectively,

In addition to this skew Markov chain, it will be convenient to construct a coupling of two copies $U_n = (X_n, Y \circ \Theta^n)$ and $U'_n = (X'_n, Y' \circ \Theta^n)$ of the skew chain such that $Y = Y'$. To construct such a coupling, we define an $E \times E \times \Omega^Y$-valued Markov process $V_n = (X_n, X'_n, Y \circ \Theta^n)$ with transition probability kernel

$$P^V(z, z', y, B \times C \times D) = P^X(z, y, B) P^X(z', y, C) I_D(\Theta y).$$

Note that the probability measure on $E \times E \times \Omega^Y$,

$$\tilde{\lambda}(A) = \int I_A(z, z', y) \mu(y, dz) \mu(y, dz') \mathbf{P}^Y(dy) = ((\mu \otimes \mu) \mathbf{P}^Y)(A),$$

is an invariant measure for the transition probability $P^V$. We will construct in the usual way a probability kernel $\mathbf{Q}_{\cdot, \cdot, \cdot} : E \times E \times \Omega^Y \times \mathcal{B}(E \times E \times \Omega^Y)^{\mathbb{Z}_+} \to [0, 1]$ such that $\mathbf{Q}_{z,z',y}$ is the law of $(V_n)_{n \geq 0}$ with $V_0 \sim \delta_{z,z',y}$. Note that under $\mathbf{Q}_{z,z',y}$, the processes $(X_n)_{n \geq 0}$ and $(X'_n)_{n \geq 0}$ are independent and their laws coincide with the law of $(X_n)_{n \geq 0}$ under $\mathbf{P}_{z,y}$ and $\mathbf{P}_{z',y}$, respectively.

Define the sequence of measurable functions

$$\beta_n(z, z', y) = \|\mathbf{P}_{z,y}(X_n \in \cdot) - \mathbf{P}_{z',y}(X_n \in \cdot)\|_{\mathrm{TV}}, \qquad n \in \mathbb{N}.$$

Note that $\beta_n$ is nonincreasing with $n$, so that $\beta(z, z', y) = \lim_{n \to \infty} \beta_n(z, z', y)$ is well defined and measurable. We wish to prove that condition 3 of Theorem 2.3 implies that $\beta(z, z', y) = 0$ $(\mu \otimes \mu) \mathbf{P}^Y$-a.e. We will do this in two steps. First, following Derriennic [16] (see also Ornstein and Sucheston [27]), we prove a zero-two law for $\beta(z, z', y)$ which asserts that either conditions 1 and 2 of Theorem 2.3 hold, or else $\beta(z, z', y)$ attains values arbitrarily close to 2. In the second step, we will show that condition 3 of Theorem 2.3 rules out the latter possibility.

PROPOSITION 2.7 (Zero-two law). *Let $\tilde{H}$ be a given set of $(\mu \otimes \mu) \mathbf{P}^Y$-full measure. Then one or the other of the following possibilities must hold true:*

1. *Condition 2 of Theorem 2.3 holds for a subset $H \subset \tilde{H}$ of $(\mu \otimes \mu)\mathbf{P}^Y$-full measure, and $\beta(z, z', y) = 0$ for all $(z, z', y) \in H$.*
2. *There is an $y \in \Omega^Y$ such that the following holds: for any $\varepsilon > 0$, there is a $(z, z', y') \in \tilde{H}$ with $y' = \Theta^n y$ for some $n \in \mathbb{N}$ and $\beta(z, z', y') > 2 - \varepsilon$.*

PROOF.  Let $H \subset \tilde{H}$ be the subset constructed through Lemma 2.6. It suffices to show that if condition 2 of Theorem 2.3 does not hold on $H$, then the second possibility in the statement of the current proposition must hold true. Indeed, if condition 2 of Theorem 2.3 does hold on $H$, then $\beta(z, z', y) = 0$ for all $(z, z', y) \in H$ by Lemma 2.5 and, thus, the first possibility holds true.

We suppose, therefore, that condition 2 of Theorem 2.3 does not hold on $H$. Then we may clearly choose a $(z, z', y) \in H$ and an $A \in \mathcal{T}^X$ such that we have either $\mathbf{P}_{z,y}(A) \neq \mathbf{P}_{z',y}(A)$ or $0 < \mathbf{P}_{z,y}(A) < 1$. Let us now define

$$Z = 2I_A - 1, \qquad g_n(\tilde{z}) = \mathbf{E}_{\tilde{z}, \Theta^n y}(Z \circ \Theta^{-n}) \qquad \text{for all } \tilde{z} \in E.$$

Using the first property of Lemma 2.2, it is not difficult to establish that

$$g_n(\tilde{z}) = \mathbf{E}_{\tilde{z}, \Theta^n y}(g_{n+k}(X_k)) \qquad \text{for all } \tilde{z} \in E, k \geq 0,$$

and that

$$g_n(X_n) = \mathbf{E}_{\tilde{z}, y}(Z | \mathcal{F}_{[0,n]}^X) \qquad \mathbf{P}_{\tilde{z}, y}\text{-a.s. for every } \tilde{z} \in E.$$

In particular, $g_n(X_n) \to Z$ $\mathbf{P}_{\tilde{z}, y}$-a.s. for every $\tilde{z} \in E$ by martingale convergence, and this implies for any $0 < \varepsilon < 2$ and $\tilde{z} \in E$ that

$$\mathbf{P}_{\tilde{z}, y}(g_n(X_n) > 1 - \varepsilon) \xrightarrow{n \to \infty} \mathbf{P}_{\tilde{z}, y}(A),$$

$$\mathbf{P}_{\tilde{z}, y}(g_n(X_n) < -1 + \varepsilon) \xrightarrow{n \to \infty} 1 - \mathbf{P}_{\tilde{z}, y}(A).$$

We now proceed as follows. Note that for any $0 < \varepsilon < 2$,

$$\mathbf{Q}_{z, z', y}(g_n(X_n) > 1 - \varepsilon/2 \text{ and } g_n(X_n') < -1 + \varepsilon/2)$$
$$= \mathbf{P}_{z, y}(g_n(X_n) > 1 - \varepsilon/2)\mathbf{P}_{z', y}(g_n(X_n) < -1 + \varepsilon/2),$$

which converges as $n \to \infty$ to $\mathbf{P}_{z, y}(A)(1 - \mathbf{P}_{z', y}(A))$, and similarly,

$$\mathbf{Q}_{z, z', y}(g_n(X_n') > 1 - \varepsilon/2 \text{ and } g_n(X_n) < -1 + \varepsilon/2)$$
$$= \mathbf{P}_{z', y}(g_n(X_n) > 1 - \varepsilon/2)\mathbf{P}_{z, y}(g_n(X_n) < -1 + \varepsilon/2),$$

which converges as $n \to \infty$ to $\mathbf{P}_{z', y}(A)(1 - \mathbf{P}_{z, y}(A))$. But as either $\mathbf{P}_{z, y}(A) \neq \mathbf{P}_{z', y}(A)$ or $0 < \mathbf{P}_{z, y}(A) < 1$, at least one of these expressions must be positive. Hence, for every $0 < \varepsilon < 2$, we can find an $n \in \mathbb{N}$ such that

$$\mathbf{Q}_{z, z', y}(|g_n(X_n) - g_n(X_n')| > 2 - \varepsilon) > 0.$$

In particular, there must then be a choice of $(\tilde{z}, \tilde{z}', \Theta^n y) \in H$ such that we have $|g_n(\tilde{z}) - g_n(\tilde{z}')| > 2 - \varepsilon$. It remains to note that, for all $k \geq 0$,

$$\beta_k(\tilde{z}, \tilde{z}', \Theta^n y) = \sup_{\|f\|_\infty \leq 1} |\mathbf{E}_{\tilde{z}, \Theta^n y}(f(X_k)) - \mathbf{E}_{\tilde{z}', \Theta^n y}(f(X_k))|$$
$$\geq |\mathbf{E}_{\tilde{z}, \Theta^n y}(g_{n+k}(X_k)) - \mathbf{E}_{\tilde{z}', \Theta^n y}(g_{n+k}(X_k))|$$
$$= |g_n(\tilde{z}) - g_n(\tilde{z}')| > 2 - \varepsilon,$$

so that $\beta(\tilde{z}, \tilde{z}', \Theta^n y) > 2 - \varepsilon$. But we can repeat this procedure for any $0 < \varepsilon < 2$, and this establishes that the second possibility of the proposition holds. $\square$

It remains to argue that condition 3 of Theorem 2.3 rules out the second possibility of the zero-two law. We will need the following lemma.

LEMMA 2.8. *The following holds for all* $(z, z', y) \in E \times E \times \Omega^Y$:

$$\beta_{n+1}(z, z', y) \leq (P^V \beta_n)(z, z', y) = \int \beta_n(\tilde{z}, \tilde{z}', \tilde{y}) P^V(z, z', y, d\tilde{z}, d\tilde{z}', d\tilde{y}).$$

*In particular,* $\beta(z, z', y) \leq (P^V \beta)(z, z', y)$.

PROOF. Choose sets $E_k^n$ as in Lemma 2.4, and define

$$\beta_\ell^n(z, z', y) = \sum_{k=1}^n |\mathbf{P}_{z,y}(X_\ell \in E_k^n) - \mathbf{P}_{z',y}(X_\ell \in E_k^n)|.$$

Then $\beta_\ell^n \nearrow \beta_\ell$ as $n \to \infty$. But $\beta_{\ell+1}^n \leq P^V \beta_\ell^n$ follows from Jensen's inequality and Lemma 2.2, so that $\beta_{\ell+1} \leq P^V \beta_\ell$ follows by monotone convergence. Letting $\ell \to \infty$, we obtain $\beta \leq P^V \beta$ by dominated convergence. $\square$

The following result now essentially completes the proof.

PROPOSITION 2.9. *Suppose that condition* 3 *of Theorem 2.3 holds. Then there is a set* $\tilde{H}$ *of* $(\mu \otimes \mu)\mathbf{P}^Y$-*full measure such that* $\beta(z, z', y) = \beta(\tilde{z}, \tilde{z}', \tilde{y}) < 2$ *for every* $(z, z', y), (\tilde{z}, \tilde{z}', \tilde{y}) \in \tilde{H}$ *with* $\tilde{y} = \Theta^n y$ *for some* $n \geq 0$.

PROOF. Denote by $\mathbf{Q}$ the law of $(V_n)_{n \geq 0}$ with initial measure $\tilde{\mu} = (\mu \otimes \mu)\mathbf{P}^Y$. By the previous lemma, $\beta(V_n)$ is a bounded submartingale under $\mathbf{Q}$ and, hence, $\{\beta(V_n)\}$ is a Cauchy sequence in $L^1(\mathbf{Q})$ by the martingale convergence theorem. But then, using the stationarity of $\mathbf{Q}$, we find that

$$\mathbf{E}_\mathbf{Q}|\beta(V_0) - \beta(V_n)| = \mathbf{E}_\mathbf{Q}|\beta(V_k) - \beta(V_{n+k})| \xrightarrow{k \to \infty} 0 \qquad \text{for all } n \in \mathbb{N}.$$

In particular, we evidently have

$$\int \mathbf{Q}_{z,z',y}(\beta(V_0) = \beta(V_n) \text{ for all } n) \tilde{\lambda}(dz, dz', dy) = 1$$

and there is consequently a set $\tilde{H}_1$ of $\tilde{\lambda}$-full measure such that

$$\mathbf{Q}_{z,z',y}(\beta(z, z', y) = \beta(V_n) \text{ for all } n) = 1 \qquad \text{for all } (z, z', y) \in \tilde{H}_1.$$

By condition 3 of Theorem 2.3, we may choose another set $\tilde{H}_2$ of $\tilde{\lambda}$-full measure such that for every $(z, \tilde{z}, y) \in \tilde{H}_2$, there is an $n \in \mathbb{N}$ such that $\mathbf{P}_{z,y}(X_n \in \cdot)$ and $\mathbf{P}_{\tilde{z},y}(X_n \in \cdot)$ are not mutually singular. Note that the latter implies that $\mathbf{P}_{z,y}(X_m \in \cdot)$ and $\mathbf{P}_{\tilde{z},y}(X_m \in \cdot)$ are not mutually singular for every $m \geq n$, as $\mathbf{P}_{z,y}(X_n \in \cdot) \perp \mathbf{P}_{\tilde{z},y}(X_n \in \cdot)$ is equivalent to $\beta_n(z, \tilde{z}, y) = 2$ and $\beta_m(z, \tilde{z}, y)$ is nonincreasing with $m$. Now define the set

$$\tilde{H}_3 = \{(z, z', \tilde{z}, \tilde{z}', y) : (z, z', y), (\tilde{z}, \tilde{z}', y) \in \tilde{H}_1, (z, \tilde{z}, y), (z', \tilde{z}', y) \in \tilde{H}_2\}.$$

Then it is easily seen that $\tilde{H}_3$ has $(\mu \otimes \mu \otimes \mu \otimes \mu)\mathbf{P}^Y$-full measure.

We claim that $\beta(z, z', y) = \beta(\tilde{z}, \tilde{z}', y)$ whenever $(z, z', \tilde{z}, \tilde{z}', y) \in \tilde{H}_3$. To see this, fix such a point, and choose $n \in \mathbb{N}$ such that $\mathbf{P}_{z,y}(X_n \in \cdot)$ and $\mathbf{P}_{\tilde{z},y}(X_n \in \cdot)$ are not mutually singular and $\mathbf{P}_{z',y}(X_n \in \cdot)$ and $\mathbf{P}_{\tilde{z}',y}(X_n \in \cdot)$ are not mutually singular. This implies, in particular, that $\mathbf{Q}_{z,z',y}(V_n \in \cdot)$ and $\mathbf{Q}_{\tilde{z},\tilde{z}',y}(V_n \in \cdot)$ are not mutually singular. But these measures are supported, respectively, on the sets

$$\Xi_1 = \{(\zeta, \zeta', \Theta^n y) : \beta(z, z', y) = \beta(\zeta, \zeta', \Theta^n y)\},$$
$$\Xi_2 = \{(\zeta, \zeta', \Theta^n y) : \beta(\tilde{z}, \tilde{z}', y) = \beta(\zeta, \zeta', \Theta^n y)\}$$

as $(z, z', y), (\tilde{z}, \tilde{z}', y) \in \tilde{H}_1$, and, as the measures are nonsingular, we must have $\Xi_1 \cap \Xi_2 \neq \varnothing$. We have therefore established that $\beta(z, z', y) = \beta(\tilde{z}, \tilde{z}', y)$.

To proceed, we define

$$\beta(y) = \int \beta(z, z', y) \mu(y, dz) \mu(y, dz').$$

We claim that $\beta(z, z', y) = \beta(y)$ $\tilde{\lambda}$-a.e. Indeed, note that

$$\int |\beta(z, z', y) - \beta(y)| \tilde{\lambda}(dz, dz', dy)$$

$$\leq \int |\beta(z, z', y) - \beta(\tilde{z}, \tilde{z}', y)| (\mu \otimes \mu \otimes \mu \otimes \mu)(y, dz, dz', d\tilde{z}, d\tilde{z}') \mathbf{P}^Y(dy)$$

by Jensen's inequality, and we may restrict the integral on the right-hand side to $\tilde{H}_3$, as this set has full measure. Thus, the left-hand side vanishes as claimed.

To complete the proof, let $\tilde{H}_4$ be a set of $\tilde{\lambda}$-full measure such that $\beta(z, z', y) = \beta(y)$ for all $(z, z', y) \in \tilde{H}_4$. Using Lemma 2.6, we can find a subset $\tilde{H}_5 \subset \tilde{H}_4$ of $\tilde{\lambda}$-full measure such that we have

$$\mathbf{Q}_{z,z',y}(V_n \in \tilde{H}_5 \text{ for all } n \geq 0) = 1 \qquad \text{for all } (z, z', y) \in \tilde{H}_5.$$

We now set $\tilde{H} = \tilde{H}_1 \cap \tilde{H}_2 \cap \tilde{H}_5$. Then evidently $\beta(z, z', y) = \beta(y) = \beta(\Theta^n y)$ for all $n \geq 0$ whenever $(z, z', y) \in \tilde{H}$, and $\beta(z, z', y) < 2$ as condition 3 of Theorem 2.3 holds for $(z, z', y) \in \tilde{H}$. The proof is easily completed.   $\square$

Let us now complete the proof of the implication $3 \Rightarrow 1, 2$ in Theorem 2.3. By the zero-two law, it suffices to show that condition 3 of Theorem 2.3 rules out the second possibility of Proposition 2.7. Assume that condition 3 of Theorem 2.3 holds, and apply the zero-two law with the set $\tilde{H}$ obtained from Proposition 2.9. If the second possibility of Proposition 2.7 holds, then there is an $y \in \Omega^Y$ and a sequence $(z_k, z_k', \Theta^{n_k} y) \in \tilde{H}$ such that $\beta(z_k, z_k', \Theta^{n_k} y) \to 2$ as $k \to \infty$. But by Proposition 2.9, $\beta(z_k, z_k', \Theta^{n_k} y) = \beta(z_1, z_1', \Theta^{n_1} y) < 2$ for all $k \geq 1$, which is a contradiction. Hence, the proof of Theorem 2.3 is complete.

## 3. Weak ergodicity of conditional Markov processes.

3.1. *The hidden Markov model.*   Throughout this paper we will operate in the same canonical setting as in Section 2. In this section, however, we will initially give a different construction of the measure $\mathbf{P}$ which makes $(X_n, Y_n)_{n \in \mathbb{Z}}$ a hidden Markov model; the *signal process* $X_n$ then plays the role of the unobserved component, while the *observation process* $Y_n$ is the observed component. Such hidden Markov structure is the usual setup in which nonlinear filtering problems are of interest. We will shortly see, however, that hidden Markov models are Markov chains in random environments in disguise, so that the results of Section 2 apply.

As before, the signal $X_n$ takes values in the Polish space $E$ and the observations $Y_n$ take values in the Polish space $F$. We proceed to construct a measure $\mathbf{P}$ on the canonical path space $(\Omega, \mathcal{F})$. The hidden Markov model consists of:

1. A probability kernel $P \colon E \times \mathcal{B}(E) \to [0, 1]$.
2. A probability measure $\pi$ on $(E, \mathcal{B}(E))$ such that

$$\int P(z, A) \pi(dz) = \pi(A) \qquad \text{for all } A \in \mathcal{B}(E).$$

3. A probability kernel $\Phi \colon E \times \mathcal{B}(F) \to [0, 1]$.

We now construct $\mathbf{P}$ as follows. For every $n \in \mathbb{N}$, we can define the probability measure $\mathbf{P}^{(n)}$ on $\mathcal{F}_{[-n,n]}$ as

$$\mathbf{P}^{(n)}(A) = \int I_A(x,y) \Phi(x(n), dy(n)) \cdots \Phi(x(-n), dy(-n))$$
$$\times P(x(n-1), dx(n)) \cdots P(x(-n), dx(-n+1)) \pi(dx(-n)).$$

Then $\mathbf{P}^{(n+1)}|_{\mathcal{F}_{[-n,n]}} = \mathbf{P}^{(n)}$, so that we can construct the probability measure

$$\mathbf{P} : \mathcal{F} \to [0,1], \qquad \mathbf{P}|_{\mathcal{F}_{[-n,n]}} = \mathbf{P}^{(n)} \qquad \text{for all } n \in \mathbb{N}$$

by the Kolmogorov extension theorem. Note that under $\mathbf{P}$, the signal $X_n$ is a stationary Markov chain with transition probability kernel $P(z, A)$ and stationary probability measure $\pi$, while, conditionally on the signal, the observations are independent at different times and $Y_n$ has law $\Phi(X_n, \cdot)$. We also remark that the joint process $(X_n, Y_n)_{n \in \mathbb{Z}}$ is easily seen to be itself a stationary Markov chain.

In addition to the probability measure $\mathbf{P}$, we introduce the probability kernel $\mathbf{P}^{\cdot} : E \times \mathcal{F}_+ \to [0,1]$ such that $\mathbf{P}^z$ is the law of $(X_n, Y_n)_{n \geq 0}$ started at $X_0 = z$ [i.e., under $\mathbf{P}^z$, the signal $(X_n)_{n \geq 0}$ is a Markov chain with transition probability kernel $P$ and initial measure $X_0 \sim \delta_z$, the observations $(Y_n)_{n \geq 0}$ are conditionally independent given the signal, and $Y_n$ has conditional law $\Phi(X_n, \cdot)$ given $\mathcal{F}_+^X$]. For any probability measure $\nu$ on $(E, \mathcal{B}(E))$, we define the probability measure

$$\mathbf{P}^\nu(A) = \int I_A(x,y) \mathbf{P}^z(dx, dy) \nu(dz) \qquad \text{for all } A \in \mathcal{F}_+.$$

Note that $\mathbf{P}^\pi$ is in fact the restriction of $\mathbf{P}$ to $\mathcal{F}_+$.

We now introduce two assumptions on the hidden Markov model which will play an important role in our main results.

ASSUMPTION 3.1 (Ergodicity). The following holds:

$$\|\mathbf{P}^z(X_n \in \cdot) - \pi\|_{\mathrm{TV}} \xrightarrow{n \to \infty} 0 \qquad \text{for } \pi\text{-a.e. } z \in E.$$

ASSUMPTION 3.2 (Nondegeneracy). There exists a probability measure $\varphi$ on $\mathcal{B}(F)$ and a strictly positive measurable function $g : E \times F \to ]0, \infty[$ such that

$$\Phi(z, A) = \int I_A(u) g(z, u) \varphi(du) \qquad \text{for all } A \in \mathcal{B}(F), z \in E.$$

We do not automatically assume in the following that either of these assumptions is in force, but we will impose them explicitly where they are needed.

3.2. *The conditional signal process.* Despite that we have constructed the measure $\mathbf{P}$ in a rather different manner, the hidden Markov model introduced in the previous subsection is in fact a disguised Markov chain in a random environment in the sense of Section 2. This is established in the following lemma.

LEMMA 3.3. *There exist probability kernels $P^X : E \times \Omega^Y \times \mathcal{B}(E) \to [0, 1]$ and $\mu : \Omega^Y \times \mathcal{B}(E) \to [0, 1]$, and a probability measure $\mathbf{P}^Y$ on $(\Omega^Y, \mathcal{F}^Y)$, such that the conditions of Section 2 are satisfied and the measure $\mathbf{P}$ constructed there coincides with the measure $\mathbf{P}$ constructed in the current section. In particular,*

$$P^X(X_n, Y \circ \Theta^n, A) = \mathbf{P}(X_{n+1} \in A | \mathcal{F}_n^X \vee \mathcal{F}^Y) \qquad \mathbf{P}\text{-}a.s.,$$

$$\mu(Y \circ \Theta^n, A) = \mathbf{P}(X_n \in A | \mathcal{F}^Y) \qquad \mathbf{P}\text{-}a.s.$$

*for every $A \in \mathcal{B}(E)$ and $n \in \mathbb{Z}$, and $\mathbf{P}^Y = \mathbf{P}|_{\mathcal{F}^Y}$.*

PROOF. Let us fix the measure $\mathbf{P}$ as defined in the current section. We will use this measure to construct $P^X$, $\mu$ and $\mathbf{P}^Y$. Subsequently, denoting by $\mathbf{P}'$ the probability measure on $\mathcal{F}$ constructed from $P^X$, $\mu$ and $\mathbf{P}^Y$ in Section 2 (called $\mathbf{P}$ there), we will show that in fact $\mathbf{P}' = \mathbf{P}$.

Set $\mathbf{P}^Y = \mathbf{P}|_{\mathcal{F}^Y}$, and let $\tilde{\mu} : \Omega^Y \times \mathcal{B}(E) \to [0, 1]$ be a regular conditional probability of the form $\mathbf{P}(X_0 \in \cdot | \mathcal{F}^Y)$. Moreover, note that

$$\mathbf{P}(X_1 \in A | \mathcal{F}_0^X \vee \mathcal{F}^Y) = \mathbf{P}(X_1 \in A | \sigma(X_0) \vee \mathcal{F}^Y) \qquad \mathbf{P}\text{-a.s.}$$

by the Markov property of $(X_n, Y_n)_{n \in \mathbb{Z}}$; indeed, due to the Markov property the $\sigma$-fields $\mathcal{F}_{[1, \infty[}$ and $\mathcal{F}_{-1}$ are conditionally independent given $\sigma(X_0, Y_0)$, so that the claim follows directly from the elementary properties of the conditional expectation. We can therefore obtain a regular conditional probability $\tilde{P}^X : E \times \Omega^Y \times \mathcal{B}(E) \to [0, 1]$ of the form $\mathbf{P}(X_1 \in \cdot | \mathcal{F}_0^X \vee \mathcal{F}^Y)$ [i.e., $\tilde{P}^X(X_0, Y, A) = \mathbf{P}(X_1 \in A | \mathcal{F}_0^X \vee \mathcal{F}^Y)$ $\mathbf{P}$-a.s. for every $A \in \mathcal{B}(E)$]. The regular conditional probabilities exist by the Polish assumption [21], Theorem 5.3.

Note that it follows trivially from the stationarity of $(X_n, Y_n)_{n \in \mathbb{Z}}$ that $\mathbf{P}^Y$ is invariant under $\Theta$. We now claim that for $\mathbf{P}^Y$-a.e. $y \in \Omega^Y$, we have

$$\int \tilde{P}^X(z, y, A) \tilde{\mu}(y, dz) = \tilde{\mu}(\Theta y, A) \qquad \text{for all } A \in \mathcal{B}(E).$$

To see this, note that as $\mathcal{B}(E)$ is countably generated, it suffices by a standard monotone class argument to prove the claim for $A$ in a countable generating algebra $\{E_n\} \subset \mathcal{B}(E)$ such that $\mathcal{B}(E) = \sigma\{E_n : n \in \mathbb{N}\}$. But note that for fixed $n \in \mathbb{N}$,

$$\int \tilde{P}^X(z, Y, E_n) \tilde{\mu}(Y, dz) = \mathbf{E}(\mathbf{P}(X_1 \in E_n | \mathcal{F}_0^X \vee \mathcal{F}^Y) | \mathcal{F}^Y) = \mathbf{P}(X_1 \in E_n | \mathcal{F}^Y),$$

while $\mathbf{P}(X_1 \in E_n | \mathcal{F}^Y) = \tilde{\mu}(Y \circ \Theta, E_n)$ follows from

$$\mathbf{E}(f(Y)\{\mathbf{P}(X_0 \in E_n | \mathcal{F}^Y) \circ \Theta\}) = \mathbf{E}(f(Y \circ \Theta^{-1})\mathbf{P}(X_0 \in E_n | \mathcal{F}^Y))$$
$$= \mathbf{E}(f(Y \circ \Theta^{-1})I_{E_n}(X_0))$$
$$= \mathbf{E}(f(Y)I_{E_n}(X_1))$$

for every bounded measurable $f : \Omega^Y \to \mathbb{R}$, where we have twice used the stationarity of $\mathbf{P}$. As we must only verify equality for a countable collection $\{E_n\}$, we can indeed find a set $H \in \mathcal{F}^Y$ of $\mathbf{P}^Y$-full measure such that

$$\int \tilde{P}^X(z, y, A)\tilde{\mu}(y, dz) = \tilde{\mu}(\Theta y, A) \qquad \text{for all } A \in \mathcal{B}(E), y \in H.$$

We now set $\mu(y, A) = \tilde{\mu}(y, A)$ and $P^X(z, y, A) = \tilde{P}^X(z, y, A)$ for all $z \in E$, $y \in H$, and $A \in \mathcal{B}(E)$, and we set $\mu(y, A) = \pi(A)$, $P^X(z, y, A) = \pi(A)$ whenever $y \notin H$. Then $\mu$ and $P^X$ are still versions of their defining regular conditional probabilities and $P^X$, $\mu$, $\mathbf{P}^Y$ satisfy the conditions of Section 2. The various identities in the statement of the lemma follow from the stationarity of $\mathbf{P}$ in the same way as we established above that $\mathbf{P}(X_1 \in E_n | \mathcal{F}^Y) = \tilde{\mu}(Y \circ \Theta, E_n)$.

It remains to show that the measure $\mathbf{P}'$ constructed from $P^X$, $\mu$, $\mathbf{P}^Y$ as in Section 2 coincides with the measure $\mathbf{P}$. It suffices to show that $\mathbf{P}'(A) = \mathbf{P}(A)$ for every $A \in \mathcal{F}_{[-n,n]}$, $n \in \mathbb{N}$. To this end, note that for $A \in \mathcal{F}_{[-n,n]}$ we evidently have

$$\mathbf{P}'(A) = \int I_A(x, y)P^X(x(n-1), \Theta^{n-1}y, dx(n)) \cdots$$
$$\times P^X(x(-n), \Theta^{-n}y, dx(-n+1))\mu(\Theta^{-n}y, dx(-n))\mathbf{P}^Y(dy)$$
$$= \mathbf{E}(\mathbf{E}(\mathbf{E}(\cdots \mathbf{E}(\mathbf{E}(I_A | \mathcal{F}_{n-1}^X \vee \mathcal{F}^Y) | \mathcal{F}_{n-2}^X \vee \mathcal{F}^Y) \cdots | \mathcal{F}_{-n}^X \vee \mathcal{F}^Y) | \mathcal{F}^Y))$$
$$= \mathbf{P}(A).$$

Thus, the proof is complete. $\quad\square$

From this point onward we will fix $P^X$, $\mu$, $\mathbf{P}^Y$ as defined in the previous lemma. In particular, this allows us to define the probability kernels $\mathbf{P}_y$ and $\mathbf{P}_{z,y}$ as in Section 2, and these are easily seen to be versions of the regular conditional probabilities $\mathbf{P}(\cdot | \mathcal{F}^Y)$ and $\mathbf{P}(\cdot | \mathcal{F}_0^X \vee \mathcal{F}^Y)$, respectively. Under $\mathbf{P}_y$, the process $(X_n)_{n \in \mathbb{Z}}$ has the law of the signal process conditioned on the observations $(Y_n)_{n \in \mathbb{Z}}$; we will refer to this process as the *conditional signal process*. The main purpose of this section is to obtain a sufficient condition for the conditional signal to be weakly ergodic, that is, for any (hence all) of the conditions of Theorem 2.3 to hold in the current setting. In Sections 4–6, we will see that this question has important consequences for the asymptotic properties of nonlinear filters.

Intuitively, it seems plausible that the weak ergodicity of the conditional signal process is inherited from the ergodicity of the (unconditional) signal process, that is, that weak ergodicity of the conditional signal follows from Assumption 3.1. The counterexample in [1] illustrates, however, that this need not be the case. The following theorem, which is the main result of this section, shows that weak ergodicity of the conditional signal follows nonetheless if we also assume nondegeneracy of the observations (Assumption 3.2).

THEOREM 3.4.   *Suppose that both Assumptions 3.1 and 3.2 are in force. Then any (hence all) of the conditions of Theorem 2.3 hold true.*

The proof of this result is contained in the following subsections.

3.3. *Weak ergodicity of the conditional signal.*   The strategy of the proof of Theorem 3.4 is to show that condition 3 of Theorem 2.3 follows from Assumptions 3.1 and 3.2. In this subsection we prove that condition 3 of Theorem 2.3 follows from Assumption 3.1 and a certain absolute continuity assumption; that the latter follows from Assumptions 3.1 and 3.2 is established in the next subsection.

LEMMA 3.5.   *Suppose Assumption 3.1 holds, and that there is a strictly positive measurable function $h : E \times \Omega^Y \times E \to ]0, \infty[$ such that for $\mu \mathbf{P}^Y$-a.e. $(z, y)$,*

$$P^X(z, y, A) = \int I_A(\tilde{z}) h(z, y, \tilde{z}) P(z, d\tilde{z}) \qquad \text{for all } A \in \mathcal{B}(E).$$

*Then condition 3 of Theorem 2.3 holds.*

PROOF.   First, we note that Assumption 3.1 implies that there is a set $H_1$ of $(\mu \otimes \mu) \mathbf{P}^Y$-full measure such that for any $(z, z', y) \in H_1$, there is an $n \in \mathbb{N}$ such that $\mathbf{P}^z(X_n \in \cdot)$ and $\mathbf{P}^{z'}(X_n \in \cdot)$ are not mutually singular. To see this, note that

$$\int \|\mathbf{P}^z(X_n \in \cdot) - \mathbf{P}^{z'}(X_n \in \cdot)\|_{\mathrm{TV}} \, \mu(y, dz) \mu(y, dz') \mathbf{P}^Y(dy)$$

$$\leq 2 \int \|\mathbf{P}^z(X_n \in \cdot) - \pi\|_{\mathrm{TV}} \, \mu(y, dz) \mathbf{P}^Y(dy)$$

$$= 2 \int \|\mathbf{P}^z(X_n \in \cdot) - \pi\|_{\mathrm{TV}} \, \pi(dz) \xrightarrow{n \to \infty} 0$$

by Assumption 3.1. But as $\|\mathbf{P}^z(X_n \in \cdot) - \mathbf{P}^{z'}(X_n \in \cdot)\|_{\mathrm{TV}}$ is nonincreasing and uniformly bounded, we find that $\|\mathbf{P}^z(X_n \in \cdot) - \mathbf{P}^{z'}(X_n \in \cdot)\|_{\mathrm{TV}} \to 0$ as $n \to \infty$ for $(\mu \otimes \mu) \mathbf{P}^Y$-a.e. $(z, z', y)$, which establishes the claim.

Now let $H_2$ be a set of $\mu \mathbf{P}^Y$-full measure such that the absolute continuity condition in the statement of the lemma holds true for all $(z, y) \in H_2$. By Lemma 2.6, there is a subset $H_3 \subset H_2$ of $\mu \mathbf{P}^Y$-full measure such that for every $(z, y) \in H_3$ we have $\mathbf{P}_{z,y}((X_n, \Theta^n y) \in H_3$ for all $n \geq 0) = 1$. It follows directly that for every $(z, y) \in H_3$, $n \in \mathbb{N}$ and $A \in \mathcal{B}(E)$, we have

$$\mathbf{P}_{z,y}(X_n \in A) = \mathbf{E}^z(h(X_0, y, X_1) \cdots h(X_{n-1}, \Theta^{n-1} y, X_n) I_A(X_n)).$$

In particular, $\mathbf{P}_{z,y}(X_n \in \cdot) \sim \mathbf{P}^z(X_n \in \cdot)$ for all $(z, y) \in H_3$ and $n \in \mathbb{N}$.

To complete the proof, define the following set:

$$H_4 = \{(z, z', y) : (z, z', y) \in H_1, (z, y), (z', y) \in H_3\}.$$

Then $H_4$ has $(\mu \otimes \mu) \mathbf{P}^Y$-full measure, and for every $(z, z', y) \in H_4$, there is an $n \in \mathbb{N}$ such that $\mathbf{P}_{z,y}(X_n \in \cdot)$ and $\mathbf{P}_{z',y}(X_n \in \cdot)$ are not mutually singular. $\square$

3.4. *Nondegeneracy.* Before we proceed, we will prove an elementary result on regular conditional probabilities. The result generalizes the trivial identity

$$\frac{\mathbf{P}(A|B, C)}{\mathbf{P}(A|C)} = \frac{\mathbf{P}(B|A, C)}{\mathbf{P}(B|C)} \qquad \text{provided } \mathbf{P}(A \cap C) > 0, \mathbf{P}(B \cap C) > 0$$

to regular conditional probabilities in Polish spaces.

LEMMA 3.6. *Let $G_1$, $G_2$ and $K$ be Polish spaces and set $\Omega = G_1 \times G_2 \times K$. We consider a probability measure $\mathbf{P}$ on $(\Omega, \mathcal{B}(\Omega))$. Denote by $\gamma_1 : \Omega \to G_1$, $\gamma_2 : \Omega \to G_2$ and $\kappa : \Omega \to K$ the coordinate projections, and let $\mathcal{G}_1$, $\mathcal{G}_2$ and $\mathcal{K}$ be the $\sigma$-fields generated by $\gamma_1$, $\gamma_2$ and $\kappa$, respectively. Choose fixed versions of the following regular conditional probabilities (which exist by the Polish assumption):*

$$\Xi_1^K(g_1, \cdot) = \mathbf{P}(\kappa \in \cdot|\mathcal{G}_1)(g_1), \qquad \Xi_{12}^K(g_1, g_2, \cdot) = \mathbf{P}(\kappa \in \cdot|\mathcal{G}_1 \vee \mathcal{G}_2)(g_1, g_2),$$

$$\Xi_1^2(g_1, \cdot) = \mathbf{P}(\gamma_2 \in \cdot|\mathcal{G}_1)(g_1), \qquad \Xi_{1K}^2(g_1, k, \cdot) = \mathbf{P}(\gamma_2 \in \cdot|\mathcal{G}_1 \vee \mathcal{K})(g_1, k),$$

*where $g_1 \in G_1$, $g_2 \in G_2$, $k \in K$. Suppose that there exists a nonnegative measurable function $h : G_1 \times G_2 \times K \to [0, \infty[$ and a set $H \subset G_1 \times G_2$ such that $\mathbf{E}(I_H(\gamma_1, \gamma_2)) = 1$ and for every $(g_1, g_2) \in H$,*

$$\Xi_{12}^K(g_1, g_2, A) = \int I_A(k) h(g_1, g_2, k) \Xi_1^K(g_1, dk) \qquad \text{for all } A \in \mathcal{K}.$$

*Then there is an $H' \subset G_1 \times K$ with $\mathbf{E}(I_{H'}(\gamma_1, \kappa)) = 1$ so that for all $(g_1, k) \in H'$,*

$$\Xi_{1K}^2(g_1, k, B) = \int I_B(g_2) h(g_1, g_2, k) \Xi_1^2(g_1, dg_2) \qquad \text{for all } B \in \mathcal{G}_2.$$

PROOF. We can evidently write (using the disintegration of measures [21], Theorem 5.4) for every $A \in \mathcal{G}_1$, $B \in \mathcal{G}_2$, and $C \in \mathcal{K}$

$$\mathbf{P}(\gamma_1 \in A, \gamma_2 \in B, \kappa \in C)$$
$$= \int I_A(g_1) I_B(g_2) \Xi_{12}^K(g_1, g_2, C) \Xi_1^2(g_1, dg_2) \Xi_1(dg_1)$$
$$= \int I_A(g_1) I_C(k) \Xi_{1K}^2(g_1, k, B) \Xi_1^K(g_1, dk) \Xi_1(dg_1),$$

where $\Xi_1$ is the law of $\gamma_1$ under $\mathbf{P}$. Therefore,

$$\int \Xi_{1K}^2(g_1, k, B) I_A(g_1) I_C(k) \Xi_1^K(g_1, dk) \Xi_1(dg_1)$$
$$= \int I_B(g_2) h(g_1, g_2, k) \Xi_1^2(g_1, dg_2) I_A(g_1) I_C(k) \Xi_1^K(g_1, dk) \Xi_1(dg_1),$$

where the exchange of integration order is permitted due to the nonnegativity of the integrand. As this holds for every $A \in \mathcal{G}_1$ and $C \in \mathcal{K}$, we obtain

$$\Xi_{1K}^2(g_1, k, B) = \int I_B(g_2) h(g_1, g_2, k) \Xi_1^2(g_1, dg_2) \qquad \text{for } \mathbf{P}\text{-a.e. } (g_1, k)$$

for every fixed $B \in \mathcal{G}_2$. But as $\mathcal{G}_2$ is countably generated, it suffices to verify that equality holds for $B$ in a countable generating algebra, and we can thus eliminate the dependence of the exceptional set on $B$. $\square$

To complete the proof of Theorem 3.4, we must show that the absolute continuity condition $P^X(z, y, \cdot) \sim P(z, \cdot)$ of Lemma 3.5 holds. Recall that $P(z, \cdot)$ is a version of the regular conditional probability $\mathbf{P}(X_1 \in \cdot | \mathcal{F}_0^X)$, while $P^X$ is a version of the regular conditional probability $\mathbf{P}(X_1 \in \cdot | \mathcal{F}_0^X \vee \mathcal{F}^Y)$. By the Markov property, however, it is immediate that we can also consider $P$ to be a version of the regular conditional probability $\mathbf{P}(X_1 \in \cdot | \sigma(X_0))$, and $P^X$ a version of the regular conditional probability $\mathbf{P}(X_1 \in \cdot | \sigma(X_0) \vee \mathcal{F}_+^Y)$. To prove absolute continuity, we will apply the previous lemma to the law of the triple $(X_0, X_1, (Y_k)_{k \geq 0})$. In particular, to establish that $P^X(z, y, \cdot) \sim P(z, \cdot)$, we may equivalently investigate whether the laws of $(Y_k)_{k \geq 0}$ under different initial conditions are equivalent.

The following result, which is of independent interest, shows that—provided the observations are nondegenerate—two initial laws of the signal give rise to equivalent laws of the observations whenever the signal forgets the initial laws. This will be used below to establish that $P^X(z, y, \cdot) \sim P(z, \cdot)$.

LEMMA 3.7. *Suppose Assumption 3.2 holds. Let $\nu, \bar{\nu}$ be probability measures such that $\|\mathbf{P}^\nu(X_n \in \cdot) - \mathbf{P}^{\bar{\nu}}(X_n \in \cdot)\|_{\mathrm{TV}} \xrightarrow{n \to \infty} 0$. Then $\mathbf{P}^\nu|_{\mathcal{F}_+^Y} \sim \mathbf{P}^{\bar{\nu}}|_{\mathcal{F}_+^Y}$.*

PROOF. We will work on the space $\Omega' = E^{\mathbb{Z}_+} \times E^{\mathbb{Z}_+} \times F^{\mathbb{Z}_+}$, where we write $X_n(x, x', y) = x(n)$, $X_n'(x, x', y) = x'(n)$, and $Y_n(x, x', y) = y(n)$.

We make use of the well-known fact [24], Theorem III.14.10 and (III.20.7), that $\|\mathbf{P}^\nu(X_n \in \cdot) - \mathbf{P}^{\bar{\nu}}(X_n \in \cdot)\|_{\mathrm{TV}} \to 0$ as $n \to \infty$ implies the existence of a successful coupling of the laws of $(X_n)_{n \geq 0}$ under $\mathbf{P}^\nu$ and $\mathbf{P}^{\bar{\nu}}$. We can thus construct a probability measure $\mathbf{Q} : \mathcal{B}(E^{\mathbb{Z}_+} \times E^{\mathbb{Z}_+}) \to [0, 1]$ such that:

1. The law of $(X_n)_{n \geq 0}$ under $\mathbf{Q}$ coincides with the law of $(X_n)_{n \geq 0}$ under $\mathbf{P}^\nu$;
2. The law of $(X_n')_{n \geq 0}$ under $\mathbf{Q}$ coincides with the law of $(X_n)_{n \geq 0}$ under $\mathbf{P}^{\bar{\nu}}$;
3. There is a finite random time $\tau$ such that a.s. $X_n = X_n'$ for all $n \geq \tau$.

In addition, we define a probability kernel $\mathbf{Q}^Y : E^{\mathbb{Z}_+} \times \mathcal{B}(F^{\mathbb{Z}_+}) \to [0, 1]$ such that $(Y_n)_{n \geq 0}$ are independent under $\mathbf{Q}^Y(x, \cdot)$ and $\mathbf{Q}^Y(x, Y_n \in \cdot) = \Phi(x(n), \cdot)$.

Now consider the following probability measures on $\Omega'$:

$$\mathbf{Q}_1(A) = \int I_A(x, x', y) \mathbf{Q}^Y(x, dy) \mathbf{Q}(dx, dx'),$$

$$\mathbf{Q}_2(A) = \int I_A(x, x', y) \mathbf{Q}^Y(x', dy) \mathbf{Q}(dx, dx').$$

It is easily seen that $\mathbf{P}^\nu|_{\mathcal{F}_+^Y} = \mathbf{Q}_1|_{\mathcal{F}_+^Y}$ and $\mathbf{P}^{\bar{\nu}}|_{\mathcal{F}_+^Y} = \mathbf{Q}_2|_{\mathcal{F}_+^Y}$. To complete the proof, it therefore suffices to show that $\mathbf{Q}_1 \sim \mathbf{Q}_2$. It is immediate, however, that

$$\frac{d\mathbf{Q}^Y(x', \cdot)}{d\mathbf{Q}^Y(x, \cdot)} = \prod_{k=0}^{N} \frac{g(x'(k), y(k))}{g(x(k), y(k))} \qquad \text{whenever } x(n) = x'(n) \text{ for all } n > N,$$

where $g(z, y)$ is the observation density defined in Assumption 3.2. Thus, evidently

$$\mathbf{Q}_1 \sim \mathbf{Q}_2 \qquad \text{with } \frac{d\mathbf{Q}_2}{d\mathbf{Q}_1} = \prod_{k=0}^{\tau} \frac{g(X_k', Y_k)}{g(X_k, Y_k)}.$$

The proof is complete. □

We can now prove the following.

LEMMA 3.8. *Suppose Assumptions 3.1 and 3.2 hold. Then there is a strictly positive measurable $h : E \times \Omega^Y \times E \to ]0, \infty[$ such that for $\mu \mathbf{P}^Y$-a.e. $(z, y)$,*

$$P^X(z, y, A) = \int I_A(\tilde{z}) h(z, y, \tilde{z}) P(z, d\tilde{z}) \qquad \text{for all } A \in \mathcal{B}(E).$$

PROOF.   By the Markov property, $P$ and $P^X$ are versions of the regular conditional probabilities $\mathbf{P}(X_1 \in \cdot | \sigma(X_0))$ and $\mathbf{P}(X_1 \in \cdot | \sigma(X_0) \vee \mathcal{F}_+^Y)$, respectively. By the Polish assumption, we can also introduce regular conditional probabilities $R : E \times \mathcal{F}_+^Y \to [0,1]$ and $R^X : E \times E \times \mathcal{F}_+^Y \to [0,1]$ of the form $\mathbf{P}((Y_k)_{k\geq 0} \in \cdot | \sigma(X_0))$ and $\mathbf{P}((Y_k)_{k\geq 0} \in \cdot | \sigma(X_0, X_1))$, respectively. Applying Lemma 3.6 to the law of the triple $(X_0, X_1, (Y_k)_{k\geq 0})$, it evidently suffices to show that there is a strictly positive measurable $h : E \times \Omega^Y \times E \to {]0, \infty[}$ such that

$$R^X(z, z', A) = \int I_A(y) h(z, y, z') R(z, dy) \qquad \text{for all } A \in \mathcal{F}_+^Y$$

for $(z, z') \in H$ with $\mathbf{P}((X_0, X_1) \in H) = 1$.

By a well-known result on kernels ([15], Section V.58) there exists a non-negative measurable function $\tilde{h} : E \times \Omega^Y \times E \to [0, \infty[$ such that, for all $z, z' \in E$,

$$R^X(z, z', A) = \int I_A(y) \tilde{h}(z, y, z') R(z, dy) + R^\perp(z, z', A) \qquad \text{for all } A \in \mathcal{F}_+^Y,$$

where the kernel $R^\perp$ is such that $R^\perp(z, z', \cdot) \perp R(z, \cdot)$ for every $z, z' \in E$. Now suppose we can establish that $R^X(z, z', \cdot) \sim R(z, \cdot)$ for $(z, z') \in H$ with $\mathbf{P}((X_0, X_1) \in H) = 1$. Then $R^\perp(z, z', \cdot) = 0$ for $(z, z') \in H$, and $\tilde{h}(z, y, z') > 0$ except on a null set. We can then set $h(z, y, z') = 1$ whenever $\tilde{h}(z, y, z') = 0$, and set $h(z, y, z') = \tilde{h}(z, y, z')$ otherwise; this gives a function $h$ with the desired properties, completing the proof. It thus remains to show that $R^X(z, z', \cdot) \sim R(z, \cdot)$ for $(z, z') \in H$ with $\mathbf{P}((X_0, X_1) \in H) = 1$.

To this end, let us introduce convenient versions of the regular conditional probabilities $R$ and $R^X$. Note that we may set

$$\int f_0(y(0)) \cdots f_n(y(n)) R^X(z, z', dy)$$

$$= \int f_0(u) \Phi(z, du) \times \mathbf{E}^{z'}(f_1(Y_0) \cdots f_n(Y_{n-1}))$$

for all bounded measurable $f_0, \ldots, f_n$ and $n < \infty$. Similarly, we may set

$$\int f_0(y(0)) \cdots f_n(y(n)) R(z, dy)$$

$$= \int f_0(u) \Phi(z, du) \times \int \mathbf{E}^{\tilde{z}}(f_1(Y_0) \cdots f_n(Y_{n-1})) P(z, d\tilde{z})$$

$$= \int f_0(u) \Phi(z, du) \times \mathbf{E}^{P(z, \cdot)}(f_1(Y_0) \cdots f_n(Y_{n-1})).$$

It thus suffices to show that

$$\mathbf{P}^{z'}|_{\mathcal{F}_+^Y} \sim \mathbf{P}^{P(z, \cdot)}|_{\mathcal{F}_+^Y} \qquad \text{for } (z, z') \in H \text{ with } \mathbf{P}((X_0, X_1) \in H) = 1.$$

By Assumption 3.2 and Lemma 3.7, it suffices to show that

$$\|\mathbf{P}^{z'}(X_n \in \cdot) - \mathbf{P}^{P(z,\cdot)}(X_n \in \cdot)\|_{\mathrm{TV}} \xrightarrow{n \to \infty} 0$$

for $(z, z') \in H$ with $\mathbf{P}((X_0, X_1) \in H) = 1$.

Now note that by Assumption 3.1, we may choose a set $H_1$ of $\pi$-full measure such that $\|\mathbf{P}^z(X_n \in \cdot) - \pi\|_{\mathrm{TV}} \to 0$ as $n \to \infty$ for all $z \in H_1$. By Lemma 2.6, there is a subset $H_2 \subset H_1$ of $\pi$-full measure such that for every $z \in H_2$ we have $\mathbf{P}^z(X_n \in H_2 \text{ for all } n \geq 0) = 1$. In particular, for $z, z' \in H_2$, we then have

$$\|\mathbf{P}^{z'}(X_n \in \cdot) - \mathbf{P}^{P(z,\cdot)}(X_n \in \cdot)\|_{\mathrm{TV}}$$
$$\leq \|\mathbf{P}^{z'}(X_n \in \cdot) - \pi\|_{\mathrm{TV}}$$
$$+ \int \|\mathbf{P}^{z''}(X_n \in \cdot) - \pi\|_{\mathrm{TV}} P(z, dz'') \xrightarrow{n \to \infty} 0.$$

But $H = H_2 \times H_2$ satisfies $\mathbf{P}((X_0, X_1) \in H) = 1$ by construction.    $\square$

Combining Lemmas 3.5 and 3.8 now completes the proof of Theorem 3.4.

## 4. Exchange of intersection and supremum of $\sigma$-fields.

As is discussed in the Introduction and in the following sections, key to the asymptotic properties of nonlinear filters are certain identities for the observation and signal $\sigma$-fields. For example, key to the proof of total variation stability (Section 5) is the identity

$$\bigcap_{n \geq 0} \mathcal{F}_+^Y \vee \mathcal{F}_{[n,\infty[}^X \overset{?}{=} \mathcal{F}_+^Y \qquad \mathbf{P}\text{-a.s.},$$

and the goal of this section is to show that such identities hold under Assumptions 3.1 and 3.2. The question can be seen as pertaining to the permissibility of the exchange of the intersection and the supremum of $\sigma$-fields; indeed, under Assumption 3.1 the tail $\sigma$-field $\mathcal{T}^X$ is $\mathbf{P}$-a.s. trivial, so that the above identity can be written as

$$\bigcap_{n \geq 0} \mathcal{F}_+^Y \vee \mathcal{F}_{[n,\infty[}^X \overset{?}{=} \mathcal{F}_+^Y \vee \bigcap_{n \geq 0} \mathcal{F}_{[n,\infty[}^X \qquad \mathbf{P}\text{-a.s.}$$

The validity of such an exchange is a notoriously delicate problem [37].

For the sake of demonstration, we begin by proving the following lemma.

LEMMA 4.1. *Suppose that any (hence all) of the conditions of Theorem 2.3 are in force. Then the following holds true:*

$$\bigcap_{n \geq 0} \mathcal{F}^Y \vee \mathcal{F}_{[n,\infty[}^X = \bigcap_{n \geq 0} \mathcal{F}^Y \vee \mathcal{F}_{-n}^X = \mathcal{F}^Y \qquad \mathbf{P}\text{-}a.s.$$

The interest of this lemma is independent of the remainder of the paper; it follows directly from Theorem 2.3, and thus serves as a simplified demonstration of the proof of the exchange of intersection and supremum property. Unfortunately, this result is not in itself of use in proving asymptotic properties of nonlinear filters, as the entire observation field $\mathcal{F}^Y$ appears in the expression rather than the positive and negative time observations $\mathcal{F}_+^Y$ and $\mathcal{F}_0^Y$. Using additional coupling and time reversal arguments, we will prove the following useful result.

THEOREM 4.2.   *Suppose that Assumptions 3.1 and 3.2 are in force. Then*

$$\bigcap_{n\geq 0} \mathcal{F}_+^Y \vee \mathcal{F}_{[n,\infty[}^X = \mathcal{F}_+^Y \quad and \quad \bigcap_{n\geq 0} \mathcal{F}_0^Y \vee \mathcal{F}_{-n}^X = \mathcal{F}_0^Y \qquad \mathbf{P}\text{-}a.s.$$

The proof of Lemma 4.1 is given in Section 4.1 below, while the proof of Theorem 4.2 is contained in Sections 4.2–4.4.

4.1. *Proof of Lemma 4.1.*   In [37], von Weizsäcker studied problems of this type in a general setting, and Lemma 4.1 can be derived from his result and Theorem 2.3. As the idea is straightforward, however, we give a direct proof here.

Let us begin by proving the assertion

$$\bigcap_{n\geq 0} \mathcal{F}^Y \vee \mathcal{F}_{[n,\infty[}^X = \mathcal{F}^Y \qquad \mathbf{P}\text{-a.s.}$$

It suffices to show that, for every $A \in \mathcal{F}$,

$$\mathbf{P}\left(A \,\Big|\, \bigcap_{n\geq 0} \mathcal{F}^Y \vee \mathcal{F}_{[n,\infty[}^X\right) = \mathbf{P}(A|\mathcal{F}^Y) \qquad \mathbf{P}\text{-a.s.}$$

As bounded random variables of the form $F(x,y) = f(x)g(y)$ are total in $L^1(\mathbf{P})$, it suffices to verify the statement for $A \in \mathcal{F}^X$ only. By the martingale convergence theorem, it is sufficient to show that, for any $A \in \mathcal{F}^X$,

$$\mathbf{P}(A|\mathcal{F}^Y \vee \mathcal{F}_{[n,\infty[}^X) \xrightarrow{n\to\infty} \mathbf{P}(A|\mathcal{F}^Y) \qquad \text{in } L^1(\mathbf{P}).$$

We now appeal to the following fact: as $\mathcal{F}_{[n,\infty[}^X$ is countably generated, we have

$$\mathbf{P}(A|\mathcal{F}^Y \vee \mathcal{F}_{[n,\infty[}^X) = \mathbf{P}_Y(A|\mathcal{F}_{[n,\infty[}^X) \qquad \mathbf{P}\text{-a.s.}$$

for any $A \in \mathcal{F}^X$, where we have used that (Lemma 2.2) $\mathbf{P}_Y(\cdot)$ is a regular conditional probability of the form $\mathbf{P}(\cdot|\mathcal{F}^Y)$; see [37], Lemma 4.II.1. But

$$\mathbf{P}_y(|\mathbf{P}_y(A|\mathcal{F}_{[n,\infty[}^X) - \mathbf{P}_y(A)|) \xrightarrow{n\to\infty} 0 \qquad \text{for } \mathbf{P}^Y\text{-a.e. } y$$

follows by martingale convergence and the following lemma.

Lemma 4.3. *Suppose that any (hence all) of the conditions of Theorem 2.3 hold. Then the tail $\sigma$-field $\mathcal{T}^X$ is $\mathbf{P}_y$-trivial for $\mathbf{P}^Y$-a.e. $y$.*

Proof. By condition 1 of Theorem 2.3, we find that

$$\int \|\mathbf{P}_{z,y}(X_n \in \cdot) - \mathbf{P}_y(X_n \in \cdot)\|_{\mathrm{TV}} \mu(y, dz) \mathbf{P}^Y(dy)$$

$$\leq \int \|\mathbf{P}_{z,y}(X_n \in \cdot) - \mathbf{P}_{z',y}(X_n \in \cdot)\|_{\mathrm{TV}} \mu(y, dz') \mu(y, dz) \mathbf{P}^Y(dy)$$

converges to zero as $n \to \infty$. But as $\|\mathbf{P}_{z,y}(X_n \in \cdot) - \mathbf{P}_y(X_n \in \cdot)\|_{\mathrm{TV}}$ is non-increasing, we find that $\|\mathbf{P}_{z,y}(X_n \in \cdot) - \mathbf{P}_y(X_n \in \cdot)\|_{\mathrm{TV}} \to 0$ as $n \to \infty$ for $\mu \mathbf{P}^Y$-a.e. $(z, y)$. Note that by the Markov property of $(X_n)_{n \geq 0}$ under $\mathbf{P}_{z,y}$,

$$\|\mathbf{P}_{z,y}(X_n \in \cdot) - \mathbf{P}_y(X_n \in \cdot)\|_{\mathrm{TV}}$$

$$= \|\mathbf{P}_{z,y}|_{\mathcal{F}^X_{[n,\infty[}} - \mathbf{P}_y|_{\mathcal{F}^X_{[n,\infty[}}\|_{\mathrm{TV}} \xrightarrow{n \to \infty} \|\mathbf{P}_{z,y}|_{\mathcal{T}^X} - \mathbf{P}_y|_{\mathcal{T}^X}\|_{\mathrm{TV}}$$

(see, e.g., [24], Section III.20). Therefore, $\mathbf{P}_{z,y}|_{\mathcal{T}^X} = \mathbf{P}_y|_{\mathcal{T}^X}$ for $\mu \mathbf{P}^Y$-a.e. $(z, y)$, and it remains to invoke condition 2 of Theorem 2.3. □

We can now easily complete the proof of $\bigcap_{n \geq 0} \mathcal{F}^Y \vee \mathcal{F}^X_{[n,\infty[} = \mathcal{F}^Y$ $\mathbf{P}$-a.s. Indeed, integrating with respect to $\mathbf{P}^Y$, we find by dominated convergence that

$$\mathbf{P}(|\mathbf{P}_Y(A|\mathcal{F}^X_{[n,\infty[}) - \mathbf{P}_Y(A)|) \xrightarrow{n \to \infty} 0$$

and the result now follows directly.

We now turn to the proof of the assertion

$$\bigcap_{n \geq 0} \mathcal{F}^Y \vee \mathcal{F}^X_{-n} = \mathcal{F}^Y \qquad \mathbf{P}\text{-a.s.}$$

As above, it suffices to show that, for every $A \in \mathcal{F}^X$,

$$\mathbf{P}(A|\mathcal{F}^Y \vee \mathcal{F}^X_{-n}) \xrightarrow{n \to \infty} \mathbf{P}(A|\mathcal{F}^Y) \qquad \text{in } L^1(\mathbf{P}).$$

In fact, it suffices to establish only that

$$\mathbf{E}(f_1(X_{k_1}) \cdots f_\ell(X_{k_\ell})|\mathcal{F}^Y \vee \mathcal{F}^X_{-n})$$

$$\xrightarrow{n \to \infty} \mathbf{E}(f_1(X_{k_1}) \cdots f_\ell(X_{k_\ell})|\mathcal{F}^Y) \qquad \text{in } L^1(\mathbf{P})$$

for all $\ell < \infty$, $k_1, \ldots, k_\ell \in \mathbb{Z}$, and bounded measurable functions $f_1, \ldots, f_\ell$, as the family of functions of the form $f_1(X_{k_1}) \cdots f_\ell(X_{k_\ell})$ is total in $L^1(\mathcal{F}^X, \mathbf{P})$. Now note that by the last property of Lemma 2.2, we can write

$$\mathbf{E}(f_1(X_{k_1}) \cdots f_\ell(X_{k_\ell})|\mathcal{F}^Y \vee \mathcal{F}^X_{-n}) = \mathbf{E}_{X_{-n}, Y \circ \Theta^{-n}}(f_1(X_{k_1+n}) \cdots f_\ell(X_{k_\ell+n})),$$

$$\mathbf{E}(f_1(X_{k_1}) \cdots f_\ell(X_{k_\ell})|\mathcal{F}^Y) = \mathbf{E}_{Y \circ \Theta^{-n}}(f_1(X_{k_1+n}) \cdots f_\ell(X_{k_\ell+n})).$$

Therefore, using the stationarity of $\mathbf{P}$, we find that

$$\mathbf{E}(|\mathbf{E}(\Lambda_0|\mathcal{F}^Y \vee \mathcal{F}_{-n}^X) - \mathbf{E}(\Lambda_0|\mathcal{F}^Y)|)$$

$$= \int |\mathbf{E}_{z,y}(\Lambda_n) - \mathbf{E}_y(\Lambda_n)|\mu(y,dz)\mathbf{P}^Y(dy)$$

$$\leq \int |\mathbf{E}_{z,y}(\Lambda_n) - \mathbf{E}_{z',y}(\Lambda_n)|\mu(y,dz)\mu(y,dz')\mathbf{P}^Y(dy),$$

where we have written $\Lambda_n = f_1(X_{k_1+n})\cdots f_\ell(X_{k_\ell+n})$ for simplicity. It follows (see, e.g., [24], Section III.20) from the first condition of Theorem 2.3 that this expression converges to zero as $n \to \infty$, and thus the claim is established.

4.2. *Time reversal.* In order to apply the theory of Markov chains in random environments, it was important to condition the signal process on all observations $\mathcal{F}^Y$. Note that the conditional probability $\mathbf{P}(X_0 \in \cdot | \mathcal{F}^Y)$ satisfies the property $\mathbf{P}(X_0 \in \cdot | \mathcal{F}^Y) \circ \Theta^n = \mathbf{P}(X_n \in \cdot | \mathcal{F}^Y)$ which was used repeatedly in Section 2; this property is not shared by the conditional probability $\mathbf{P}(X_0 \in \cdot | \mathcal{F}_+^Y)$. An unfortunate consequence is that we obtain the triviality of $\mathcal{T}^X$ under the regular conditional probability $\mathbf{P}(\cdot | \mathcal{F}^Y)$, which leads to Lemma 4.1, rather than the triviality of $\mathcal{T}^X$ under $\mathbf{P}(\cdot | \mathcal{F}_+^Y)$, which would give (the first part of) Theorem 4.2.

To prove Theorem 4.2, we must therefore eliminate the dependence of our results to date on the past observations. As we will see in the following subsections, this can be done provided that the signal is not only ergodic forward in time (as is guaranteed by Assumption 3.1) but also after time reversal; in essence, we aim to establish that the remote past of the signal does not depend on the present. In this subsection, we will show that this property in fact already follows from Assumption 3.1, so that no additional assumptions need to be imposed.

In the following we will extend the definition of $\mathbf{P}^z$ to negative times, that is, $\mathbf{P}^z$ is a version of the regular conditional probability $\mathbf{P}(\cdot | X_0)$. Note that the time reversed signal $\tilde{X}_n = X_{-n}$ is again a Markov chain under $\mathbf{P}$ and $\mathbf{P}^z$ with stationary measure $\pi$. The goal of this subsection is to prove the following result.

PROPOSITION 4.4. *Suppose that Assumption 3.1 holds. Then*

$$\|\mathbf{P}^z(X_{-n} \in \cdot) - \pi\|_{\mathrm{TV}} \xrightarrow{n\to\infty} 0 \qquad \text{for } \pi\text{-a.e. } z \in E.$$

We will need the following lemma on regular conditional probabilities.

LEMMA 4.5. *Let $G$ be a Polish space. Denote by $\gamma_1 : G \times G \to G$ and $\gamma_2 : G \times G \to G$ the coordinate projections and by $\mathcal{G}_1$ and $\mathcal{G}_2$ the $\sigma$-fields* *generated by $\gamma_1$ and $\gamma_2$, respectively. Consider a probability measure $\pi$ on $(G, \mathcal{B}(G))$, and a probability measure $\mathbf{P}$ on $(G \times G, \mathcal{B}(G \times G))$ such that the laws of $\gamma_1$ and $\gamma_2$ under $\mathbf{P}$ both equal $\pi$. Denote by $P_1 : G \times \mathcal{B}(G) \to [0,1]$ and $P_2 : G \times \mathcal{B}(G) \to [0,1]$ the regular conditional probabilities of the form $\mathbf{P}(\gamma_1 \in \cdot | \mathcal{G}_2)$ and $\mathbf{P}(\gamma_2 \in \cdot | \mathcal{G}_1)$, respectively, and consider their Lebesgue decompositions*

$$\mathbf{P}(A \times B) = \int I_A(dz) I_B(dz') p(z, z') \pi(dz) \pi(dz') + \mathbf{P}^\perp(A \times B),$$

$$P_1(z', A) = \int I_A(z) p_1(z, z') \pi(dz) + P_1^\perp(z', A),$$

$$P_2(z, B) = \int I_B(z') p_2(z, z') \pi(dz') + P_2^\perp(z, B),$$

*where $\mathbf{P}^\perp \perp \pi \otimes \pi$, $P_1^\perp(z', \cdot) \perp \pi$ and $P_2^\perp(z, \cdot) \perp \pi$, and $p, p_1, p_2 : G \times G \to [0, \infty[$ are measurable. Then $p(z, z') = p_1(z, z') = p_2(z, z')$ for $\pi \otimes \pi$-a.e. $(z, z')$.*

PROOF.  The existence of regular conditional probabilities follows from the Polish assumption, while the existence of measurable $p_1, p_2$ follows from [15], Section V.58. It also follows from [15], Sections V.56–58, that there exist $S_1, S_2 \in \mathcal{B}(G \times G)$ such that $(\pi \otimes \pi)(S_1) = (\pi \otimes \pi)(S_2) = 1$ and for $\pi$-a.e. $z, z'$,

$$\int I_{S_1}(z, z') P_1^\perp(z', dz) = 0, \qquad \int I_{S_2}(z, z') P_2^\perp(z, dz') = 0.$$

Now note that, by the disintegration of measures, we have for all $A, B \in \mathcal{B}(G)$

$$\mathbf{P}(A \times B) = \int I_B(z') P_1(z', A) \pi(dz') = \int I_A(z) P_2(z, B) \pi(dz).$$

Now substitute in the Lebesgue decompositions of $P_1$ and $P_2$, and note that

$$\int I_{S_1}(z, z') P_1^\perp(z', dz) \pi(dz') = \int I_{S_2}(z, z') P_2^\perp(z, dz') \pi(dz) = 0.$$

Therefore, $P_1^\perp \pi \perp \pi \otimes \pi$ and $P_2^\perp \pi \perp \pi \otimes \pi$. But by the uniqueness of the Lebesgue decomposition of $\mathbf{P}$, this implies that

$$\int I_A(dz) I_B(dz') p(z, z') \pi(dz) \pi(dz')$$

$$= \int I_A(dz) I_B(dz') p_1(z, z') \pi(dz) \pi(dz')$$

$$= \int I_A(dz) I_B(dz') p_2(z, z') \pi(dz) \pi(dz')$$

for all $A, B \in \mathcal{B}(G)$, from which the result follows.  $\square$

We can now prove Proposition 4.4.

PROOF OF PROPOSITION 4.4.   Denote by $f_n(z, z')$ the density in the Lebesgue decomposition of $\mathbf{P}^z(X_n \in \cdot)$ with respect to $\pi$. Then by Assumption 3.1,

$$\int |f_n(z, z') - 1| \pi(dz) \pi(dz') \xrightarrow{n \to \infty} 0.$$

In particular, there is a subsequence $n_k \nearrow \infty$ such that

$$\int |f_{n_k}(z, z') - 1| \pi(dz) \xrightarrow{k \to \infty} 0 \qquad \text{for } \pi\text{-a.e. } z'.$$

But by the previous lemma and by stationarity, $f_n(z, z')$ is also the density in the Lebesgue decomposition of $\mathbf{P}^{z'}(X_{-n} \in \cdot)$ with respect to $\pi$. It follows that $\|\mathbf{P}^{z'}(X_{-n_k} \in \cdot) - \pi\|_{\mathrm{TV}} \to 0$ as $k \to \infty$ for $\pi$-a.e. $z'$. But $\tilde{X}_n = X_{-n}$ is again Markov, so $\|\mathbf{P}^{z'}(X_{-n} \in \cdot) - \pi\|_{\mathrm{TV}}$ is nonincreasing and the result follows.   □

4.3. *Equivalence of the initial measures.*   Let us begin by fixing a version $\mu^+ : \Omega^Y \times \mathcal{B}(E) \to [0, 1]$ of the regular conditional probability $\mathbf{P}(X_0 \in \cdot | \mathcal{F}_+^Y)$. We can then define a probability kernel $\mathbf{P}^+ : \Omega^Y \times \mathcal{F}_+^X \to [0, 1]$ by setting

$$\mathbf{P}_y^+(A) = \int \mathbf{P}_{z,y}(A) \mu^+(y, dz) \qquad \text{for all } A \in \mathcal{F}_+^X, y \in \Omega^Y.$$

It is not difficult to see that $\mathbf{P}_y^+$ is a version of the regular conditional probability $\mathbf{P}(\cdot | \mathcal{F}_+^Y)$; indeed, it suffices to note that by the Markov property $\mathbf{P}_{z,y}$ is a version of the regular conditional probability $\mathbf{P}(\cdot | \sigma(X_0) \vee \mathcal{F}_+^Y)$. We also recall that

$$\mathbf{P}_y(A) = \int \mathbf{P}_{z,y}(A) \mu(y, dz) \qquad \text{for all } A \in \mathcal{F}_+^X, y \in \Omega^Y$$

is a version of the regular conditional probability $\mathbf{P}(\cdot | \mathcal{F}^Y)$.

Theorem 2.3 establishes that the tail $\sigma$-field $\mathcal{T}^X$ is $\mathbf{P}_y$-a.s. trivial for $\mathbf{P}^Y$-a.e. $y$ (Lemma 4.3). To demonstrate the first part of Theorem 4.2 along the lines of the proof of Lemma 4.1, however, we would have to show that $\mathcal{T}^X$ is $\mathbf{P}_y^+$-a.s. trivial for $\mathbf{P}^Y$-a.e. $y$. The latter would follow from the former if we could show that $\mathbf{P}_y^+ \sim \mathbf{P}_y$ for $\mathbf{P}^Y$-a.e. $y$, and it evidently suffices to show that $\mu^+(y, \cdot) \sim \mu(y, \cdot)$ for $\mathbf{P}^Y$-a.e. $y$. The purpose of this subsection is to prove that this is indeed the case under Assumptions 3.1 and 3.2. In fact, we will prove the following stronger statement: $\mu^+(y, \cdot) \sim \pi$ and $\mu(y, \cdot) \sim \pi$ for $\mathbf{P}^Y$-a.e. $y$.

The easy part of the proof is contained in the following lemma.

Lemma 4.6. *Suppose Assumptions 3.1 and 3.2 hold. Then there is a strictly positive measurable $k^+ : \Omega^Y \times E \to ]0, \infty[$ such that, for $\mathbf{P}^Y$-a.e. $y \in \Omega^Y$,*

$$\mu^+(y, A) = \int I_A(\tilde{z}) k^+(y, \tilde{z}) \pi(d\tilde{z}) \qquad \text{for all } A \in \mathcal{B}(E).$$

Proof. By Lemma 3.6, it suffices to show that there exists a strictly positive measurable $k^+ : \Omega^Y \times E \to ]0, \infty[$ such that, for $\pi$-a.e. $z \in E$,

$$\mathbf{P}^z(B) = \int I_B(y) k^+(y, z) \mathbf{P}(dy) \qquad \text{for all } B \in \mathcal{F}_+^Y.$$

But this follows immediately from Lemma 3.7 and Assumptions 3.1 and 3.2. $\square$

It remains to prove the corresponding result for $\mu$. Though we will proceed along the same lines, the proof is complicated by the fact that Lemma 3.7 only establishes equivalence for observations at positive times $\mathcal{F}_+^Y$ and not on the entire time interval $\mathcal{F}^Y$. We therefore set out to extend Lemma 3.7 to $\mathcal{F}^Y$.

Lemma 4.7. *Under Assumptions 3.1 and 3.2, $\mathbf{P}^z|_{\mathcal{F}^Y} \sim \mathbf{P}|_{\mathcal{F}^Y}$ for $\pi$-a.e. $z$.*

Proof. By the Markov property of the signal process, $\mathcal{F}_{[n,\infty[}^X$ and $\mathcal{F}_{-n}^X$ are independent under $\mathbf{P}^z$. We can therefore estimate as follows:

$$\|\mathbf{P}^z|_{\mathcal{F}_{-n}^X \vee \mathcal{F}_{[n,\infty[}^X} - \mathbf{P}^{z'}|_{\mathcal{F}_{-n}^X \vee \mathcal{F}_{[n,\infty[}^X}\|_{\mathrm{TV}}$$

$$= \|\mathbf{P}^z|_{\mathcal{F}_{-n}^X} \otimes \mathbf{P}^z|_{\mathcal{F}_{[n,\infty[}^X} - \mathbf{P}^{z'}|_{\mathcal{F}_{-n}^X} \otimes \mathbf{P}^{z'}|_{\mathcal{F}_{[n,\infty[}^X}\|_{\mathrm{TV}}$$

$$\leq \|\mathbf{P}^z|_{\mathcal{F}_{-n}^X} - \mathbf{P}^{z'}|_{\mathcal{F}_{-n}^X}\|_{\mathrm{TV}} + \|\mathbf{P}^z|_{\mathcal{F}_{[n,\infty[}^X} - \mathbf{P}^{z'}|_{\mathcal{F}_{[n,\infty[}^X}\|_{\mathrm{TV}}$$

$$= \|\mathbf{P}^z(X_{-n} \in \cdot) - \mathbf{P}^{z'}(X_{-n} \in \cdot)\|_{\mathrm{TV}} + \|\mathbf{P}^z(X_n \in \cdot) - \mathbf{P}^{z'}(X_n \in \cdot)\|_{\mathrm{TV}}$$

$$\leq \|\mathbf{P}^z(X_{-n} \in \cdot) - \pi\|_{\mathrm{TV}} + \|\mathbf{P}^{z'}(X_{-n} \in \cdot) - \pi\|_{\mathrm{TV}}$$

$$\quad + \|\mathbf{P}^z(X_n \in \cdot) - \pi\|_{\mathrm{TV}} + \|\mathbf{P}^{z'}(X_n \in \cdot) - \pi\|_{\mathrm{TV}}.$$

Here we have used the Markov property of $X_n$ and $\tilde{X}_n = X_{-n}$, and the elementary identity $\|\mu_1 \otimes \nu_1 - \mu_2 \otimes \nu_2\|_{\mathrm{TV}} \leq \|\mu_1 - \mu_2\|_{\mathrm{TV}} + \|\nu_1 - \nu_2\|_{\mathrm{TV}}$. By Assumption 3.1 and Proposition 4.4, we now find that

$$\|\mathbf{P}^z|_{\mathcal{F}_{-n}^X \vee \mathcal{F}_{[n,\infty[}^X} - \mathbf{P}^{z'}|_{\mathcal{F}_{-n}^X \vee \mathcal{F}_{[n,\infty[}^X}\|_{\mathrm{TV}} \xrightarrow{n \to \infty} 0 \qquad \text{for } \pi \otimes \pi\text{-a.e. } (z, z').$$

But then we have

$$\|\mathbf{P}^z|_{\mathcal{F}^X_{-n} \vee \mathcal{F}^X_{[n,\infty[}} - \mathbf{P}|_{\mathcal{F}^X_{-n} \vee \mathcal{F}^X_{[n,\infty[}}\|_{\mathrm{TV}}$$

$$\leq \int \|\mathbf{P}^z|_{\mathcal{F}^X_{-n} \vee \mathcal{F}^X_{[n,\infty[}} - \mathbf{P}^{z'}|_{\mathcal{F}^X_{-n} \vee \mathcal{F}^X_{[n,\infty[}}\|_{\mathrm{TV}} \pi(dz') \xrightarrow{n \to \infty} 0 \qquad \text{for } \pi\text{-a.e. } z.$$

In particular, $\mathbf{P}$ and $\mathbf{P}^z$ agree on the remote $\sigma$-field for $\pi$-a.e. $z$:

$$\mathbf{P}^z|_{\mathcal{R}^X} = \mathbf{P}|_{\mathcal{R}^X} \qquad \text{for } \pi\text{-a.e. } z, \qquad \mathcal{R}^X = \bigcap_{n \geq 0} \mathcal{F}^X_{-n} \vee \mathcal{F}^X_{[n,\infty[}.$$

From this point onward, we fix an arbitrary $z$ such that $\mathbf{P}^z|_{\mathcal{R}^X} = \mathbf{P}|_{\mathcal{R}^X}$. To complete the proof, it suffices to show that this implies $\mathbf{P}^z|_{\mathcal{F}^Y} \sim \mathbf{P}|_{\mathcal{F}^Y}$.

To proceed, we note that the remote $\sigma$-field $\mathcal{R}^X$ coincides with the tail $\sigma$-field of the one-sided sequence $(X_{-n}, X_n)_{n \geq 0}$. We can therefore apply the maximal coupling theorem [24], Theorem III.14.10, to this sequence. In particular, we find that we can construct a probability measure $\mathbf{Q} \colon \mathcal{B}(E^{\mathbb{Z}} \times E^{\mathbb{Z}}) \to [0,1]$ such that:

1. The law of $(X_n)_{n \in \mathbb{Z}}$ under $\mathbf{Q}$ coincides with the law of $(X_n)_{n \in \mathbb{Z}}$ under $\mathbf{P}^z$;
2. The law of $(X'_n)_{n \in \mathbb{Z}}$ under $\mathbf{Q}$ coincides with the law of $(X_n)_{n \in \mathbb{Z}}$ under $\mathbf{P}$;
3. There is a random time $0 \leq \tau < \infty$ such that a.s. $X_n = X'_n$ for all $|n| \geq \tau$.

Here $X_n$ and $X'_n$ are the canonical coordinate processes on $E^{\mathbb{Z}} \times E^{\mathbb{Z}}$. The remainder of the proof now proceeds exactly as the proof of Lemma 3.7. □

We can now prove the equivalence of $\mu(y, \cdot)$ and $\pi$.

Lemma 4.8. *Suppose Assumptions 3.1 and 3.2 hold. Then there is a strictly positive measurable $k \colon \Omega^Y \times E \to ]0, \infty[$ such that, for $\mathbf{P}^Y$-a.e. $y \in \Omega^Y$,*

$$\mu(y, A) = \int I_A(\tilde{z}) k(y, \tilde{z}) \pi(d\tilde{z}) \qquad \text{for all } A \in \mathcal{B}(E).$$

Proof. By Lemma 3.6, it suffices to show that there exists a strictly positive measurable $k \colon \Omega^Y \times E \to ]0, \infty[$ such that, for $\pi$-a.e. $z \in E$,

$$\mathbf{P}^z(B) = \int I_B(y) k(y, z) \mathbf{P}(dy) \qquad \text{for all } B \in \mathcal{F}^Y.$$

But this follows immediately from Lemma 4.7 and Assumptions 3.1 and 3.2. □

The following corollary follows directly.

COROLLARY 4.9. *Suppose that Assumptions 3.1 and 3.2 hold true. Then*

$$\mathbf{P}_y^+|_{\mathcal{F}_+^X} \sim \mathbf{P}_y|_{\mathcal{F}_+^X} \qquad \text{for } \mathbf{P}^Y\text{-a.e. } y \in \Omega^Y.$$

*In particular, $\mathbf{P}_y^+|_{\mathcal{T}^X} \sim \mathbf{P}_y|_{\mathcal{T}^X}$ for $\mathbf{P}^Y$-a.e. $y \in \Omega^Y$.*

4.4. *Proof of Theorem 4.2.* We begin by proving the first assertion

$$\bigcap_{n \geq 0} \mathcal{F}_+^Y \vee \mathcal{F}_{[n,\infty[}^X = \mathcal{F}_+^Y \qquad \mathbf{P}\text{-a.s.}$$

This would follow exactly as in the proof of the first part of Lemma 4.1 if we could show that $\mathcal{T}^X$ is $\mathbf{P}_y^+$-a.s. trivial for $\mathbf{P}^Y$-a.e. $y$. But this follows directly from Lemma 4.3 and Corollary 4.9, so the claim is established.

We now turn to the second assertion

$$\bigcap_{n \geq 0} \mathcal{F}_0^Y \vee \mathcal{F}_{-n}^X = \mathcal{F}_0^Y \qquad \mathbf{P}\text{-a.s.}$$

Note that this assertion is precisely equivalent to the first assertion of the theorem after time reversal. But by Proposition 4.4, the reversed Markov chain $\tilde{X}_n = X_{-n}$ satisfies Assumption 3.1 whenever the forward chain $X_n$ does, and Assumption 3.2 is invariant under time reversal. Thus, it suffices to apply the first part of the theorem to the hidden Markov model obtained by replacing the forward transition kernel $P(z, \cdot)$ by the backward transition kernel $\mathbf{P}^z(X_{-1} \in \cdot)$. This completes the proof.

## 5. Total variation stability of the nonlinear filter.

Let us begin with a brief reminder of elementary filtering theory. The purpose of nonlinear filtering is to compute conditional probabilities of the form $\mathbf{P}^\mu(X_n \in \cdot | \mathcal{F}_{[0,n]}^Y)$. We will choose fixed versions of these regular conditional probabilities according to the following well-known lemma, whose proof we provide for future reference.

LEMMA 5.1. *Suppose that Assumption 3.2 holds. For every probability measure $\mu$ on $\mathcal{B}(E)$, we define a sequence of probability kernels $\Pi_n^\mu : \Omega^Y \times \mathcal{B}(E) \to [0, 1]$ $(n \geq 0)$ through the following recursion:*

$$\Pi_n^\mu(y, A) = \frac{\int I_A(z) g(z, y(n)) P(z', dz) \Pi_{n-1}^\mu(y, dz')}{\int g(z, y(n)) P(z', dz) \Pi_{n-1}^\mu(y, dz')},$$

$$\Pi_0^\mu(y, A) = \frac{\int I_A(z) g(z, y(0)) \mu(dz)}{\int g(z, y(0)) \mu(dz)},$$

*where $g$ is the observation density defined in Assumption 3.2. Then $\Pi_n^\mu$ is a version of the regular conditional probability $\mathbf{P}^\mu(X_n \in \cdot | \mathcal{F}_{[0,n]}^Y)$ for every $n \geq 0$.*

PROOF. Writing out the recursion, we find that

$$\Pi_n^\mu(y, A) = \frac{\mathbf{E}^\mu(g(X_0, y(0)) \cdots g(X_n, y(n)) I_A(X_n))}{\mathbf{E}^\mu(g(X_0, y(0)) \cdots g(X_n, y(n)))}.$$

But note that, by construction,

$$g(X_0, Y_0) \cdots g(X_n, Y_n) = \frac{d\mathbf{P}^\mu|_{\mathcal{F}_{[0,n]}}}{d(\mathbf{P}^\mu|_{\mathcal{F}_{[0,n]}^X} \otimes \varphi^{\otimes n})},$$

so that by the Bayes formula $\Pi_n^\mu(Y, A) = \mathbf{P}^\mu(X_n \in A | \mathcal{F}_{[0,n]}^Y)$ $\mathbf{P}^\mu$-a.s. $\quad\square$

The filter stability problem can now be phrased as follows: under which conditions does the filter $\Pi_n^\mu$ become independent of $\mu$ for large $n$? The main goal of this section is to give a precise answer to this question under Assumptions 3.1 and 3.2. To this end, we will prove the following theorem.

THEOREM 5.2. *Suppose that Assumptions 3.1 and 3.2 hold. Then*

$$\|\Pi_n^\mu - \Pi_n^\pi\|_{\mathrm{TV}} \xrightarrow{n\to\infty} 0 \qquad \mathbf{P}^\mu\text{-a.s.} \quad iff \quad \|\mathbf{P}^\mu(X_n \in \cdot) - \pi\|_{\mathrm{TV}} \xrightarrow{n\to\infty} 0.$$

The following corollaries are essentially immediate.

COROLLARY 5.3. *Suppose that Assumptions 3.1 and 3.2 hold, and call the probability measure $\mu$ stable if $\|\Pi_n^\mu - \Pi_n^\pi\|_{\mathrm{TV}} \to 0$ $\mathbf{P}^\mu$-a.s. as $n \to \infty$. Then $\mu$ is stable whenever $\mu \ll \pi$, and $\delta_z$ is stable for $\pi$-a.e. $z \in E$. Moreover, stability holds for all $\mu$ if and only if the signal process is Harris recurrent and aperiodic.*

PROOF. The first two statements follow directly from Assumption 3.1, while the last statement follows from [30], Proposition 3.6, and the fact that, by assumption, the signal possesses a finite invariant measure $\pi$. $\quad\square$

COROLLARY 5.4. *Suppose that Assumptions 3.1 and 3.2 hold true. If we have $\|\mathbf{P}^\mu(X_n \in \cdot) - \pi\|_{\mathrm{TV}} \to 0$, then $\|\Pi_n^\mu - \Pi_n^\pi\|_{\mathrm{TV}} \to 0$ $\mathbf{P}$-a.s. In particular, if*

$$\|\mathbf{P}^\mu(X_n \in \cdot) - \pi\|_{\mathrm{TV}} \xrightarrow{n\to\infty} 0, \qquad \|\mathbf{P}^\nu(X_n \in \cdot) - \pi\|_{\mathrm{TV}} \xrightarrow{n\to\infty} 0,$$

*we find that $\|\Pi_n^\mu - \Pi_n^\nu\|_{\mathrm{TV}} \to 0$ $\mathbf{P}$-a.s., $\mathbf{P}^\mu$-a.s. and $\mathbf{P}^\nu$-a.s.*

PROOF. Apply Lemma 3.7 and the triangle inequality. $\quad\square$

COROLLARY 5.5. *Suppose that Assumption 3.2 holds and that the signal is Harris recurrent and aperiodic. Then $\|\Pi_n^\mu - \Pi_n^\nu\|_{\mathrm{TV}} \to 0$ $\mathbf{P}^\gamma$-a.s. for all $\mu, \nu, \gamma$.*

PROOF.   It is well known that for Harris recurrent aperiodic Markov chains which possess a finite invariant measure $\pi$, we have $\|\mathbf{P}^\mu(X_n \in \cdot) - \pi\|_{\mathrm{TV}} \to 0$ as $n \to \infty$ for every probability measure $\mu$ [29], Theorem 6.2.8. Therefore, Assumption 3.1 follows, and it remains to apply the previous corollary and Lemma 3.7.  □

The remainder of this section is devoted to the proof of Theorem 5.2.

5.1. *Proof of Theorem 5.2: the case $\mu \ll \pi$.*   We begin by proving stability of probability measures $\mu$ that are absolutely continuous with respect to the stationary measure $\pi$. Note that by Assumption 3.1 we have $\|\mathbf{P}^\mu(X_n \in \cdot) - \pi\|_{\mathrm{TV}} \to 0$ as $n \to \infty$ for any $\mu \ll \pi$. We will also need the following result.

LEMMA 5.6.   *Suppose that Assumption 3.2 holds true and that $\mu \ll \pi$. Then we have $\Pi_n^\mu(y, \cdot) \ll \Pi_n^\pi(y, \cdot)$ for every $y \in \Omega^Y$, where*

$$\frac{d\Pi_n^\mu}{d\Pi_n^\pi}(Y, X_n) = \frac{\mathbf{E}((d\mu/d\pi)(X_0)|\mathcal{F}_+^Y \vee \mathcal{F}_{[n,\infty[}^X)}{\mathbf{E}((d\mu/d\pi)(X_0)|\mathcal{F}_{[0,n]}^Y)} \qquad \mathbf{P}\text{-}a.s.$$

PROOF.   That $\Pi_n^\mu(y, \cdot) \ll \Pi_n^\pi(y, \cdot)$ for every $y \in \Omega^Y$ can be read off directly from the expression in the proof of Lemma 5.1. Now note that

$$\frac{d\mathbf{P}^\mu}{d\mathbf{P}}\bigg|_{\mathcal{F}_{[0,\infty[}} = \frac{d\mu}{d\pi}(X_0), \qquad \frac{d\mathbf{P}^\mu}{d\mathbf{P}}\bigg|_{\mathcal{F}_{[0,n]}^Y} = \mathbf{E}\bigg(\frac{d\mu}{d\pi}(X_0)\bigg|\mathcal{F}_{[0,n]}^Y\bigg).$$

Moreover, it follows easily from Assumption 3.2 that

$$\mathbf{P}^\mu|_{\mathcal{F}_{[0,n]}^Y} \sim \mathbf{P}|_{\mathcal{F}_{[0,n]}^Y} \qquad \text{for every } n \in \mathbb{N}.$$

Therefore, the conditional expectations $\mathbf{P}^\mu(X_n \in A|\mathcal{F}_{[0,n]}^Y)$ are $\mathbf{P}$-a.s. uniquely defined and $\mathbf{E}(\frac{d\mu}{d\pi}(X_0)|\mathcal{F}_{[0,n]}^Y) > 0$ $\mathbf{P}$-a.s. We obtain by the Bayes formula

$$\mathbf{P}^\mu(X_n \in A|\mathcal{F}_{[0,n]}^Y)$$

$$= \frac{\mathbf{E}(I_A(X_n)(d\mu/d\pi)(X_0)|\mathcal{F}_{[0,n]}^Y)}{\mathbf{E}((d\mu/d\pi)(X_0)|\mathcal{F}_{[0,n]}^Y)}$$

$$= \frac{\mathbf{E}(I_A(X_n)\mathbf{E}((d\mu/d\pi)(X_0)|\sigma(X_n) \vee \mathcal{F}_{[0,n]}^Y)|\mathcal{F}_{[0,n]}^Y)}{\mathbf{E}((d\mu/d\pi)(X_0)|\mathcal{F}_{[0,n]}^Y)} \qquad \mathbf{P}\text{-a.s.}$$

Choose a measurable $\Lambda_n \colon \Omega^Y \times E \to [0, \infty[$ such that

$$\frac{\mathbf{E}((d\mu/d\pi)(X_0)|\sigma(X_n) \vee \mathcal{F}_{[0,n]}^Y)}{\mathbf{E}((d\mu/d\pi)(X_0)|\mathcal{F}_{[0,n]}^Y)} = \Lambda_n(Y, X_n) \qquad \mathbf{P}\text{-a.s.}$$

Then evidently for every $A \in \mathcal{B}(E)$

$$\Pi_n^\mu(Y, A) = \int I_A(z) \Lambda_n(Y, z) \Pi_n^\pi(Y, dz) \qquad \mathbf{P}\text{-a.s.}$$

But as $\mathcal{B}(E)$ is countably generated, it suffices by a monotone class argument to restrict to $A$ in a countable generating algebra, and we can therefore eliminate the dependence of the $\mathbf{P}$-null set on $A$. It remains to note that

$$\mathbf{E}\left(\frac{d\mu}{d\pi}(X_0)\Big|\mathcal{F}_+^Y \vee \mathcal{F}_{[n,\infty[}^X\right) = \mathbf{E}\left(\frac{d\mu}{d\pi}(X_0)\Big|\sigma(X_n) \vee \mathcal{F}_{[0,n]}^Y\right) \qquad \mathbf{P}\text{-a.s.}$$

by the Markov property, and the proof is complete. $\square$

We immediately obtain the following corollary.

COROLLARY 5.7. *Suppose Assumption 3.2 holds and $\mu \ll \pi$. Then $\mathbf{P}$-a.s.*

$$\|\Pi_n^\mu - \Pi_n^\pi\|_{\mathrm{TV}}$$
$$= \frac{\mathbf{E}(|\mathbf{E}((d\mu/d\pi)(X_0)|\mathcal{F}_+^Y \vee \mathcal{F}_{[n,\infty[}^X) - \mathbf{E}((d\mu/d\pi)(X_0)|\mathcal{F}_{[0,n]}^Y)||\mathcal{F}_{[0,n]}^Y)}{\mathbf{E}((d\mu/d\pi)(X_0)|\mathcal{F}_{[0,n]}^Y)}.$$

PROOF. This follows directly from the identity

$$\|\Pi_n^\mu(y, \cdot) - \Pi_n^\pi(y, \cdot)\|_{\mathrm{TV}} = \int \left|\frac{d\Pi_n^\mu}{d\Pi_n^\pi}(y, z) - 1\right| \Pi_n^\pi(y, dz)$$

and the previous lemma. $\square$

We can now complete the proof of Theorem 5.2 for the case $\mu \ll \pi$.

LEMMA 5.8. *Suppose Assumptions 3.1 and 3.2 hold and $\mu \ll \pi$. Then*

$$\|\Pi_n^\mu - \Pi_n^\pi\|_{\mathrm{TV}} \xrightarrow{n\to\infty} 0 \qquad \mathbf{P}\text{-a.s.}$$

*and therefore also $\mathbf{P}^\mu$-a.s. as $\mathbf{P}^\mu \ll \mathbf{P}$.*

PROOF. We aim to establish the $\mathbf{P}$-a.s. limit of the expression in Corollary 5.7. Note that the denominator satisfies

$$\mathbf{E}\left(\frac{d\mu}{d\pi}(X_0)\Big|\mathcal{F}_{[0,n]}^Y\right) \xrightarrow{n\to\infty} \mathbf{E}\left(\frac{d\mu}{d\pi}(X_0)\Big|\mathcal{F}_+^Y\right) = \frac{d\mathbf{P}^\mu}{d\mathbf{P}}\Big|_{\mathcal{F}_+^Y} \qquad \mathbf{P}\text{-a.s.}$$

by martingale convergence. Moreover, $\mathbf{P}^\mu|_{\mathcal{F}_+^Y} \sim \mathbf{P}|_{\mathcal{F}_+^Y}$ by Lemma 3.7 and Assumptions 3.1 and 3.2. Therefore, the $\mathbf{P}$-a.s. limit of the denominator is $\mathbf{P}$-a.s. strictly positive. It remains to establish convergence of the numerator.

To this end, note that for any $k \in \mathbb{N}$ we have $\mathbf{P}$-a.s.

$$\left| \mathbf{E}\left( \frac{d\mu}{d\pi}(X_0) \Big| \mathcal{F}_+^Y \vee \mathcal{F}_{[n,\infty[}^X \right) - \mathbf{E}\left( \frac{d\mu}{d\pi}(X_0) \Big| \mathcal{F}_{[0,n]}^Y \right) \right|$$

$$\leq \left| \mathbf{E}\left( \frac{d\mu}{d\pi}(X_0) I_{(d\mu/d\pi) \leq k}(X_0) \Big| \mathcal{F}_+^Y \vee \mathcal{F}_{[n,\infty[}^X \right) \right.$$
$$\left. - \mathbf{E}\left( \frac{d\mu}{d\pi}(X_0) I_{(d\mu/d\pi) \leq k}(X_0) \Big| \mathcal{F}_{[0,n]}^Y \right) \right|$$
$$+ \left| \mathbf{E}\left( \frac{d\mu}{d\pi}(X_0) I_{(d\mu/d\pi) > k}(X_0) \Big| \mathcal{F}_+^Y \vee \mathcal{F}_{[n,\infty[}^X \right) \right.$$
$$\left. - \mathbf{E}\left( \frac{d\mu}{d\pi}(X_0) I_{(d\mu/d\pi) > k}(X_0) \Big| \mathcal{F}_{[0,n]}^Y \right) \right|$$
$$\leq \left| \mathbf{E}\left( \frac{d\mu}{d\pi}(X_0) I_{(d\mu/d\pi) \leq k}(X_0) \Big| \mathcal{F}_+^Y \vee \mathcal{F}_{[n,\infty[}^X \right) \right.$$
$$\left. - \mathbf{E}\left( \frac{d\mu}{d\pi}(X_0) I_{(d\mu/d\pi) \leq k}(X_0) \Big| \mathcal{F}_{[0,n]}^Y \right) \right|$$
$$+ \mathbf{E}\left( \frac{d\mu}{d\pi}(X_0) I_{(d\mu/d\pi) > k}(X_0) \Big| \mathcal{F}_+^Y \vee \mathcal{F}_{[n,\infty[}^X \right)$$
$$+ \mathbf{E}\left( \frac{d\mu}{d\pi}(X_0) I_{(d\mu/d\pi) > k}(X_0) \Big| \mathcal{F}_{[0,n]}^Y \right).$$

In particular, setting for notational convenience

$$M_n^k = \left| \mathbf{E}\left( \frac{d\mu}{d\pi}(X_0) I_{(d\mu/d\pi) \leq k}(X_0) \Big| \mathcal{F}_+^Y \vee \mathcal{F}_{[n,\infty[}^X \right) \right.$$
$$\left. - \mathbf{E}\left( \frac{d\mu}{d\pi}(X_0) I_{(d\mu/d\pi) \leq k}(X_0) \Big| \mathcal{F}_{[0,n]}^Y \right) \right|,$$

we find that the numerator $R_n$ satisfies

$$R_n = \mathbf{E}\left( \left| \mathbf{E}\left( \frac{d\mu}{d\pi}(X_0) \Big| \mathcal{F}_+^Y \vee \mathcal{F}_{[n,\infty[}^X \right) - \mathbf{E}\left( \frac{d\mu}{d\pi}(X_0) \Big| \mathcal{F}_{[0,n]}^Y \right) \right| \Big| \mathcal{F}_{[0,n]}^Y \right)$$
$$\leq \mathbf{E}(M_n^k | \mathcal{F}_{[0,n]}^Y) + 2\mathbf{E}\left( \frac{d\mu}{d\pi}(X_0) I_{(d\mu/d\pi) > k}(X_0) \Big| \mathcal{F}_{[0,n]}^Y \right).$$

But $\mathbf{E}(M_n^k | \mathcal{F}_{[0,n]}^Y) \to 0$ $\mathbf{P}$-a.s. as $n \to \infty$ by Hunt's lemma [15], Theorem V.45, as $M_n^k \leq k$ for all $n$ and $M_n^k \to 0$ $\mathbf{P}$-a.s. as $n \to \infty$ by martingale convergence and Theorem 4.2. Moreover, by martingale convergence and dominated convergence,

$$\limsup_{k \to \infty} \limsup_{n \to \infty} \mathbf{E}\left( \frac{d\mu}{d\pi}(X_0) I_{(d\mu/d\pi) > k}(X_0) \Big| \mathcal{F}_{[0,n]}^Y \right) = 0 \qquad \mathbf{P}\text{-a.s.}$$

Therefore, the numerator converges to zero $\mathbf{P}$-a.s., and the proof is complete. $\square$

REMARK 5.9.   Along the same lines, one can prove the following result. Suppose that Assumptions 3.1 and 3.2 hold and that the relative entropy of $\mu$ with respect to $\pi$ is finite, that is, $D(\mu\|\pi) < \infty$. Then $D(\Pi_n^\mu\|\Pi_n^\pi) \to 0$ $\mathbf{P}$-a.s. as $n \to \infty$. We refer to [9] for further details on the role of relative entropy in filter stability.

5.2. *Proof of Theorem 5.2: the general case.*   We now devote our attention to the case where $\mu$ is not necessarily absolutely continuous with respect to $\pi$. Let us begin by proving the only if part of the theorem.

LEMMA 5.10.   *Suppose that Assumptions 3.1 and 3.2 hold and that*
$$\limsup_{n\to\infty} \|\mathbf{P}^\mu(X_n \in \cdot) - \pi\|_{\mathrm{TV}} > 0.$$
*Then we must have*
$$\mathbf{P}^\mu\Big(\limsup_{n\to\infty} \|\Pi_n^\mu - \Pi_n^\pi\|_{\mathrm{TV}} = 0\Big) < 1.$$

PROOF.   Let $\mathbf{P}^\mu(X_n \in \cdot) = \mu_n + \mu_n^\perp$ be the Lebesgue decomposition of $\mathbf{P}^\mu(X_n \in \cdot)$ with respect to $\pi$. In particular, $\mu_n \ll \pi$ and $\mu_n^\perp \perp \pi$, and there exists a set $S_n$ such that $\pi(S_n) = 0$ and $\mu_n^\perp(S_n^c) = 0$. We claim that
$$\limsup_{n\to\infty} \|\mathbf{P}^\mu(X_n \in \cdot) - \pi\|_{\mathrm{TV}} > 0 \quad \implies \quad \limsup_{n\to\infty} \mathbf{P}^\mu(X_n \in S_n) > 0.$$
Indeed, by [28], Theorem 7.2, Assumption 3.1 and $\mathbf{P}^\mu(X_n \in S_n) \to 0$ as $n \to \infty$ would imply that $\|\mathbf{P}^\mu(X_n \in \cdot) - \pi\|_{\mathrm{TV}} \to 0$ as $n \to \infty$, which is a contradiction.

Now note that it is easily established, using the expression in the proof of Lemma 5.1, that Assumption 3.2 implies $\Pi_n^\pi(y, \cdot) \sim \pi$ for every $y \in \Omega^Y$. Therefore, $\Pi_n^\pi(y, S_n) = 0$ for all $y \in \Omega^Y$, and we can estimate as follows:
$$\Pi_n^\pi(y, S_n) = |\Pi_n^\mu(y, S_n) - \Pi_n^\pi(y, S_n)| \le \|\Pi_n^\mu(y, \cdot) - \Pi_n^\pi(y, \cdot)\|_{\mathrm{TV}}.$$
In particular, we find that
$$\mathbf{P}^\mu(X_n \in S_n) = \mathbf{E}^\mu(\Pi_n^\mu(Y, S_n)) \le \mathbf{E}^\mu(\|\Pi_n^\mu(Y, \cdot) - \Pi_n^\pi(Y, \cdot)\|_{\mathrm{TV}})$$
and we must therefore have
$$\limsup_{n\to\infty} \mathbf{E}^\mu(\|\Pi_n^\mu(Y, \cdot) - \Pi_n^\pi(Y, \cdot)\|_{\mathrm{TV}}) > 0.$$
The proof is easily completed.   $\square$

It remains to prove the converse assertion. The idea is to reduce the general case to the case $\mu \ll \pi$. To this end, we will need the following lemma.

LEMMA 5.11.   *Suppose that Assumption 3.2 holds. Let $\mu$ and $\rho$ be probability measures, and let $\mu = \mu_{ac} + \mu_s$ be the Lebesgue decomposition of $\mu$ with respect to $\rho$ (i.e., $\mu_{ac} \ll \rho$ and $\mu_s \perp \rho$). Choose $S$ so that $\rho(S) = 1$ and $\mu_s(S) = 0$. Then*

$$\Pi_n^\mu(Y, A) = \mathbf{P}^\mu(X_0 \in S | \mathcal{F}_{[0,n]}^Y) \Pi_n^\nu(Y, A) + \mathbf{P}^\mu(X_0 \notin S | \mathcal{F}_{[0,n]}^Y) \Pi_n^{\nu^\perp}(Y, A)$$

$\mathbf{P}^\mu$-*a.s. for every $A \in \mathcal{B}(E)$, where we have written $\nu = \mu_{ac}/\mu_{ac}(E)$ and $\nu^\perp = \mu_s/\mu_s(E)$. In particular, we obtain $\mathbf{P}^\mu$-a.s. the estimate*

$$\|\Pi_n^\mu(Y, \cdot) - \Pi_n^\rho(Y, \cdot)\|_{\mathrm{TV}} \leq \|\Pi_n^\nu(Y, \cdot) - \Pi_n^\rho(Y, \cdot)\|_{\mathrm{TV}} + 2\mathbf{P}^\mu(X_0 \notin S | \mathcal{F}_{[0,n]}^Y).$$

PROOF.   Note that $d\nu/d\mu = I_S/\mu_{ac}(E)$. By the Bayes formula, we thus have

$$\mathbf{E}^\mu(I_S(X_0) I_A(X_n) | \mathcal{F}_{[0,n]}^Y) = \mathbf{E}^\mu(I_S(X_0) | \mathcal{F}_{[0,n]}^Y) \mathbf{E}^\nu(I_A(X_n) | \mathcal{F}_{[0,n]}^Y) \qquad \mathbf{P}^\mu\text{-a.s.}$$

Similarly, as $d\nu^\perp/d\mu = I_{S^c}/\mu_s(E)$, we find that

$$\mathbf{E}^\mu(I_{S^c}(X_0) I_A(X_n) | \mathcal{F}_{[0,n]}^Y) = \mathbf{E}^\mu(I_{S^c}(X_0) | \mathcal{F}_{[0,n]}^Y) \mathbf{E}^{\nu^\perp}(I_A(X_n) | \mathcal{F}_{[0,n]}^Y) \qquad \mathbf{P}^\mu\text{-a.s.}$$

The first claim now follows by summing these expressions. To prove the second assertion, let $I_k = \{E_1^k, \ldots, E_k^k\}$ be an increasing sequence of partitions of $E$ as in the proof of Lemma 2.4. Then we can estimate

$$\sum_{\ell=1}^k |\Pi_n^\mu(Y, E_\ell^k) - \Pi_n^\rho(Y, E_\ell^k)|$$

$$\leq \mathbf{P}^\mu(X_0 \in S | \mathcal{F}_{[0,n]}^Y) \sum_{\ell=1}^k |\Pi_n^\nu(Y, E_\ell^k) - \Pi_n^\rho(Y, E_\ell^k)|$$

$$+ \mathbf{P}^\mu(X_0 \notin S | \mathcal{F}_{[0,n]}^Y) \sum_{\ell=1}^k (\Pi_n^{\nu^\perp}(Y, E_\ell^k) + \Pi_n^\rho(Y, E_\ell^k))$$

$$\leq \sum_{\ell=1}^k |\Pi_n^\nu(Y, E_\ell^k) - \Pi_n^\rho(Y, E_\ell^k)| + 2\mathbf{P}^\mu(X_0 \notin S | \mathcal{F}_{[0,n]}^Y) \qquad \mathbf{P}^\mu\text{-a.s.}$$

It remains to take the limit as $k \to \infty$.   □

Note that in this result $\nu \ll \rho$ by construction. In particular, presuming that Assumptions 3.1 and 3.2 hold true and that $\|\mathbf{P}^\mu(X_n \in \cdot) - \pi\|_{\mathrm{TV}} \to 0$, and substituting $\pi$ for $\rho$, it is not difficult to establish using Lemmas 5.8 and 3.7 that

$$\limsup_{n \to \infty} \|\Pi_n^\mu(Y, \cdot) - \Pi_n^\pi(Y, \cdot)\|_{\mathrm{TV}} \leq 2\mathbf{P}^\mu(X_0 \notin S | \mathcal{F}_+^Y) \qquad \mathbf{P}^\mu\text{-a.s.}$$

We can therefore eliminate the absolutely continuous part of the initial measure $\mu$ using the stability for the case $\mu \ll \pi$ (Lemma 5.8). However, the singular part leaves the residual quantity $\mathbf{P}^\mu(X_0 \notin S | \mathcal{F}_+^Y)$, and it remains to eliminate this term. To resolve this problem, we will exploit the recursive property of the filter. Together with Lemma 5.10, the following result completes the proof of Theorem 5.2.

LEMMA 5.12. *Suppose that Assumptions 3.1 and 3.2 hold and that*

$$\limsup_{n\to\infty} \|\mathbf{P}^\mu(X_n \in \cdot) - \pi\|_{\mathrm{TV}} = 0.$$

*Then we must have*

$$\limsup_{n\to\infty} \|\Pi_n^\mu - \Pi_n^\pi\|_{\mathrm{TV}} = 0 \qquad \mathbf{P}^\mu\text{-}a.s.$$

PROOF. Define the following probability kernels:

$$\Upsilon_0^\mu(y, A) = \mu(A), \qquad \Upsilon_n^\mu(y, A) = \int I_A(z) P(z', dz) \Pi_{n-1}^\mu(y, dz').$$

Then by Lemma 5.1, the filter satisfies the recursive property

$$\Pi_{n+k}^\mu(y, A) = \Pi_k^{\Upsilon_n^\mu(y, \cdot)}(\Theta^n y, A) \qquad \text{for all } k, n \in \mathbb{Z}_+, y \in \Omega^Y, A \in \mathcal{B}(E).$$

In particular, we can write

$$\limsup_{k\to\infty} \|\Pi_k^\mu(y, \cdot) - \Pi_k^\pi(y, \cdot)\|_{\mathrm{TV}}$$

$$= \limsup_{k\to\infty} \|\Pi_k^{\Upsilon_n^\mu(y, \cdot)}(\Theta^n y, \cdot) - \Pi_k^{\Upsilon_n^\pi(y, \cdot)}(\Theta^n y, \cdot)\|_{\mathrm{TV}} \qquad \text{for all } n \in \mathbb{Z}_+.$$

But from routine manipulations, it follows that, for any $B \in \mathcal{F}_{[0,\infty[}$,

$$\mathbf{E}^\mu(I_B \circ \Theta^n | \mathcal{F}_{[0,n-1]}^Y) = \mathbf{P}^{\Upsilon_n^\mu(Y, \cdot)}(B) \qquad \mathbf{P}^\mu\text{-a.s.}$$

Therefore,

$$\mathbf{E}^\mu\Big(\limsup_{k\to\infty} \|\Pi_k^\mu(Y, \cdot) - \Pi_k^\pi(Y, \cdot)\|_{\mathrm{TV}} | \mathcal{F}_{[0,n-1]}^Y\Big)$$

$$= \mathbf{E}^\mu\Big(\limsup_{k\to\infty} \|\Pi_k^{\Upsilon_n^\mu(Y, \cdot)}(Y \circ \Theta^n, \cdot) - \Pi_k^{\Upsilon_n^\pi(Y, \cdot)}(Y \circ \Theta^n, \cdot)\|_{\mathrm{TV}} | \mathcal{F}_{[0,n-1]}^Y\Big)$$

$$= \mathbf{E}^\mu\Big(\limsup_{k\to\infty} \|\Pi_k^{\Upsilon_n^\mu(y, \cdot)}(Y \circ \Theta^n, \cdot) - \Pi_k^{\Upsilon_n^\pi(y, \cdot)}(Y \circ \Theta^n, \cdot)\|_{\mathrm{TV}} | \mathcal{F}_{[0,n-1]}^Y\Big)\Big|_{y=Y}$$

$$= \mathbf{E}^{\Upsilon_n^\mu(y, \cdot)}\Big(\limsup_{k\to\infty} \|\Pi_k^{\Upsilon_n^\mu(y, \cdot)}(Y, \cdot) - \Pi_k^{\Upsilon_n^\pi(y, \cdot)}(Y, \cdot)\|_{\mathrm{TV}}\Big)\Big|_{y=Y} \qquad \mathbf{P}^\mu\text{-a.s.,}$$

where we have used that $\Upsilon_n^\mu(Y, \cdot)$ is $\mathcal{F}_{[0,n-1]}^Y$-measurable.

For the time being, let us fix a $y \in \Omega^Y$. Note that it is easily established, using the expression in the proof of Lemma 5.1, that $\Upsilon_n^\rho(y,\cdot) \sim \mathbf{P}^\rho(X_n \in \cdot)$ for every $\rho, n, y$. Denote by $\mathbf{P}^\mu(X_n \in \cdot) = \mu_n + \mu_n^\perp$ the Lebesgue decomposition of $\mathbf{P}^\mu(X_n \in \cdot)$ with respect to $\pi$ (i.e., $\mu_n \ll \pi$ and $\mu_n^\perp \perp \pi$), and choose $S_n$ such that $\pi(S_n) = 1$ and $\mu_n^\perp(S_n) = 0$. Then clearly $\Upsilon_n^\mu(y,\cdot) = \nu_n(y,\cdot) + \nu_n^\perp(y,\cdot)$ with

$$\nu_n(y, A) = \Upsilon_n^\mu(y, A \cap S_n), \qquad \nu_n^\perp(y, A) = \Upsilon_n^\mu(y, A \cap S_n^c)$$

is the Lebesgue decomposition of $\Upsilon_n^\mu(y,\cdot)$ with respect to $\Upsilon_n^\pi(y,\cdot)$ [i.e., $\nu_n(y,\cdot) \ll \Upsilon_n^\pi(y,\cdot)$ and $\nu_n^\perp(y,\cdot) \perp \Upsilon_n^\pi(y,\cdot)$]. By Lemma 5.11, we find that

$$\|\Pi_k^{\Upsilon_n^\mu(y,\cdot)}(Y,\cdot) - \Pi_k^{\Upsilon_n^\pi(y,\cdot)}(Y,\cdot)\|_{\mathrm{TV}}$$
$$\leq \|\Pi_k^{\nu_n(y,\cdot)}(Y,\cdot) - \Pi_k^{\Upsilon_n^\pi(y,\cdot)}(Y,\cdot)\|_{\mathrm{TV}} + 2\mathbf{P}^{\Upsilon_n^\mu(y,\cdot)}(X_0 \notin S_n | \mathcal{F}_{[0,k]}^Y)$$
$$\leq \|\Pi_k^{\nu_n(y,\cdot)}(Y,\cdot) - \Pi_k^\pi(Y,\cdot)\|_{\mathrm{TV}} + \|\Pi_k^{\Upsilon_n^\pi(y,\cdot)}(Y,\cdot) - \Pi_k^\pi(Y,\cdot)\|_{\mathrm{TV}}$$
$$+ 2\mathbf{P}^{\Upsilon_n^\mu(y,\cdot)}(X_0 \notin S_n | \mathcal{F}_{[0,k]}^Y) \qquad \mathbf{P}^{\Upsilon_n^\mu(y,\cdot)}\text{-a.s.}$$

But $\nu_n(y,\cdot) \ll \pi$ and $\Upsilon_n^\pi(y,\cdot) \sim \pi$, so by Lemma 5.8 the first two terms on the right converge to zero as $k \to \infty$ $\mathbf{P}$-a.s. We claim that this convergence also holds $\mathbf{P}^{\Upsilon_n^\mu(y,\cdot)}$-a.s. Indeed, recall that $\Upsilon_n^\mu(y,\cdot) \sim \mathbf{P}^\mu(X_n \in \cdot) := \rho_n$, so that the claim is established if we can show that $\mathbf{P}^{\rho_n}|_{\mathcal{F}_+^Y} \sim \mathbf{P}|_{\mathcal{F}_+^Y}$. But $\|\mathbf{P}^{\rho_n}(X_k \in \cdot) - \pi\|_{\mathrm{TV}} = \|\mathbf{P}^\mu(X_{n+k} \in \cdot) - \pi\|_{\mathrm{TV}} \to 0$, so the claim follows from Lemma 3.7.

We have now established that, for every $y \in \Omega^Y$,

$$\mathbf{E}^{\Upsilon_n^\mu(y,\cdot)}\Big(\limsup_{k\to\infty}\|\Pi_k^{\Upsilon_n^\mu(y,\cdot)}(Y,\cdot) - \Pi_k^{\Upsilon_n^\pi(y,\cdot)}(Y,\cdot)\|_{\mathrm{TV}}\Big) \leq 2\mathbf{P}^{\Upsilon_n^\mu(y,\cdot)}(X_0 \notin S_n).$$

In particular, this implies that $\mathbf{P}^\mu$-a.s.

$$\mathbf{E}^\mu\Big(\limsup_{k\to\infty}\|\Pi_k^\mu(Y,\cdot) - \Pi_k^\pi(Y,\cdot)\|_{\mathrm{TV}} | \mathcal{F}_{[0,n-1]}^Y\Big) \leq 2\mathbf{P}^\mu(X_n \notin S_n | \mathcal{F}_{[0,n-1]}^Y)$$

and, therefore, we have for all $n \in \mathbb{N}$

$$\mathbf{E}^\mu\Big(\limsup_{k\to\infty}\|\Pi_k^\mu(Y,\cdot) - \Pi_k^\pi(Y,\cdot)\|_{\mathrm{TV}}\Big) \leq 2\mathbf{P}^\mu(X_n \notin S_n) = 2\mu_n^\perp(E).$$

But by the assumption that $\|\mathbf{P}^\mu(X_n \in \cdot) - \pi\|_{\mathrm{TV}} \to 0$, we must have $\mu_n^\perp(E) \to 0$ as $n \to \infty$. Thus, the proof is complete. $\square$

**6. Continuous time.**

6.1. *The hidden Markov model in continuous time.* Up to this point we have exclusively dealt with Markov chains and hidden Markov models in discrete time. In this section, we will prove analogous results for continuous time filtering models by reducing them to the discrete time setting. First, however, we must introduce the class of continuous time models in which we will be interested.

We consider an $\tilde{E}$-valued signal $(\xi_t)_{t \in \mathbb{R}}$ and an $\tilde{F}$-valued observation $(\eta_t)_{t \in \mathbb{R}}$, where $\tilde{E}$ is a Polish space and $\tilde{F}$ is a Polish topological vector space. We will realize these processes on the canonical path space $\tilde{\Omega} = \tilde{\Omega}^\xi \times \tilde{\Omega}^\eta$, where $\tilde{\Omega}^\xi = D(\mathbb{R}; \tilde{E})$ and $\tilde{\Omega}^\eta = D(\mathbb{R}; \tilde{F})$ are, respectively, the Skorohod spaces of $\tilde{E}$- and $\tilde{F}$-valued càdlàg paths. Denote by $\tilde{\mathscr{F}}$ the Borel $\sigma$-field on $\tilde{\Omega}$, and we introduce the natural filtrations $\tilde{\mathscr{F}}_t^\xi$, $\tilde{\mathscr{F}}_t^\eta$, $\tilde{\mathscr{F}}_t$ in complete analogy with the discrete time case:

$$\tilde{\mathscr{F}}_t^\xi = \sigma\{\xi_s : s \leq t\}, \qquad \tilde{\mathscr{F}}_t^\eta = \sigma\{\eta_s : s \leq t\}, \qquad \tilde{\mathscr{F}}_t = \tilde{\mathscr{F}}_t^\xi \vee \tilde{\mathscr{F}}_t^\eta.$$

Moreover, we define for intervals $[s,t]$ $(s \leq t)$ the $\sigma$-fields

$$\tilde{\mathscr{F}}_{[s,t]}^\xi = \sigma\{\xi_r : r \in [s,t]\}, \qquad \tilde{\mathscr{F}}_{[s,t]}^\eta = \sigma\{\eta_r - \eta_s : r \in [s,t]\}$$

and we set $\tilde{\mathscr{F}}_{[s,t]} = \tilde{\mathscr{F}}_{[s,t]}^\xi \vee \tilde{\mathscr{F}}_{[s,t]}^\eta$. Finally, we define

$$\tilde{\mathscr{F}}^\xi = \bigvee_{t \geq 0} \tilde{\mathscr{F}}_t^\xi, \qquad \tilde{\mathscr{F}}_+^\xi = \bigvee_{t \geq 0} \tilde{\mathscr{F}}_{[0,t]}^\xi, \qquad \tilde{\mathscr{F}}^\eta = \bigvee_{t \geq 0} \tilde{\mathscr{F}}_t^\eta, \qquad \tilde{\mathscr{F}}_+^\eta = \bigvee_{t \geq 0} \tilde{\mathscr{F}}_{[0,t]}^\eta.$$

The canonical shift is defined as $\tilde{\Theta}^s(\xi, \eta)(t) = (\xi(s+t), \eta(s+t) - \eta(s))$.

The continuous time hidden Markov model now consists of the following:

1. A probability kernel $\tilde{\mathbf{Q}}^\cdot : \tilde{E} \times \tilde{\mathscr{F}}_+^\xi \to [0,1]$ such that, for every $A \in \mathcal{B}(\tilde{E})$,

   $$\tilde{\mathbf{Q}}^z(\xi_t \in A | \tilde{\mathscr{F}}_s) = \tilde{\mathbf{Q}}^{\xi_s}(\xi_{t-s} \in A) \qquad \tilde{\mathbf{Q}}^z\text{-a.s. for all } z \in \tilde{E}, t \geq s \geq 0,$$

   and such that $\tilde{\mathbf{Q}}^z(\xi_0 = z) = 1$ for all $z \in \tilde{E}$.

2. A probability measure $\tilde{\pi}$ such that

   $$\int \tilde{\mathbf{Q}}^z(\xi_t \in A)\tilde{\pi}(dz) = \tilde{\pi}(A) \qquad \text{for all } A \in \mathcal{B}(\tilde{E}), t \geq 0.$$

3. A probability kernel $\tilde{\Phi} : \tilde{\Omega}^\xi \times \tilde{\mathscr{F}}^\eta \to [0,1]$ such that $(\eta_t)_{t \in \mathbb{R}}$ has independent increments with respect to $\tilde{\Phi}(\xi, \cdot)$ for every $\xi \in \tilde{\Omega}^\xi$ and such that

   $$\int I_A(\tilde{\Theta}^s \eta)\tilde{\Phi}(\xi, d\eta) = \tilde{\Phi}(\tilde{\Theta}^s \xi, A) \qquad \text{for all } \xi \in \tilde{\Omega}^\xi, A \in \tilde{\mathscr{F}}^\eta, s \in \mathbb{R}.$$

We assume, moreover, that $\tilde{\Phi}(\xi, A)$ is $\tilde{\mathscr{F}}_{[s,t]}^\xi$-measurable for every $A \in \tilde{\mathscr{F}}_{[s,t]}^\eta$.

For any probability measure $\mu$ on $\mathcal{B}(\tilde{E})$, we define

$$\tilde{\mathbf{Q}}^\mu(A) = \int \tilde{\mathbf{Q}}^z(A)\mu(dz) \qquad \text{for all } A \in \tilde{\mathcal{F}}^\xi_+.$$

Then under $\tilde{\mathbf{Q}}^\mu$, the signal $(\xi_t)_{t\geq 0}$ is a time homogeneous Markov process with initial measure $\xi_0 \sim \mu$. In particular, under $\tilde{\mathbf{Q}}^{\tilde{\pi}}$ the signal is a stationary Markov process with stationary measure $\tilde{\pi}$. We can therefore extend the measure $\tilde{\mathbf{Q}}^{\tilde{\pi}}$ to two-sided time $\tilde{\mathcal{F}}^\xi$ in the usual fashion, and we denote this extended measure as $\tilde{\mathbf{Q}}$. In particular, under $\tilde{\mathbf{Q}}$ the entire signal $(\xi_t)_{t\in\mathbb{R}}$ is a stationary Markov process with stationary measure $\tilde{\pi}$. We now define the probability measure $\tilde{\mathbf{P}}$ on $\tilde{\mathcal{F}}$ as

$$\tilde{\mathbf{P}}(A) = \int I_A(\xi,\eta)\tilde{\Phi}(\xi,d\eta)\tilde{\mathbf{Q}}(d\xi) \qquad \text{for all } A \in \tilde{\mathcal{F}}$$

and we similarly define the measures $\tilde{\mathbf{P}}^\mu$ on $\tilde{\mathcal{F}}^\xi_+ \vee \tilde{\mathcal{F}}^\eta_+$ as

$$\tilde{\mathbf{P}}^\mu(A) = \int I_A(\xi,\eta)\tilde{\Phi}(\xi,d\eta)\tilde{\mathbf{Q}}^\mu(d\xi) \qquad \text{for all } A \in \tilde{\mathcal{F}}^\xi_+ \vee \tilde{\mathcal{F}}^\eta_+.$$

Then $\tilde{\mathbf{P}}^\mu$ defines the hidden Markov model with initial measure $\mu$, while $\tilde{\mathbf{P}}$ defines the stationary hidden Markov model. Note that the stationary measure $\tilde{\mathbf{P}}$ is invariant under the canonical shift $\tilde{\Theta}^s$ by construction.

We now introduce the continuous time counterparts of Assumptions 3.1 and 3.2.

ASSUMPTION 6.1 (Ergodicity).   The following holds:

$$\|\tilde{\mathbf{Q}}^z(\xi_t \in \cdot) - \tilde{\pi}\|_{\mathrm{TV}} \xrightarrow{t\to\infty} 0 \qquad \text{for } \tilde{\pi}\text{-a.e. } z \in \tilde{E}.$$

ASSUMPTION 6.2 (Nondegeneracy).   There exists a probability measure $\tilde{\varphi}$ on $\tilde{\mathcal{F}}^\eta$ and a family $(\tilde{\Sigma}_{s,t})_{s\leq t}$ of strictly positive random variables such that

$$\tilde{\Phi}(\xi,A) = \int I_A(\eta)\tilde{\Sigma}_{s,t}(\xi,\eta)\tilde{\varphi}(d\eta) \qquad \text{for all } A \in \tilde{\mathcal{F}}^\eta_{[s,t]}, \xi \in \tilde{\Omega}^\xi, s \leq t,$$

and such that $\tilde{\Sigma}_{s,t}$ is $\tilde{\mathcal{F}}_{[s,t]}$-measurable for every $s \leq t$.

Our guiding example in which a kernel $\tilde{\Phi}$ can be constructed that satisfies all the required conditions is the ubiquitous filtering model with white noise observations. Though none of our results rely specifically on this model, let us take a moment to show that it does indeed fit within our general framework.

EXAMPLE 6.3 (White noise observations).    Set $\tilde{F} = \mathbb{R}^d$ for some $d < \infty$, and let $\tilde{\varphi}$ be the probability measure which makes $(\eta_t)_{t \in \mathbb{R}}$ a two-sided $d$-dimensional Wiener process. Such a probability measure is easily constructed; indeed, let $\mathbf{W}$ be the canonical Wiener measure on $C([0, \infty[; \mathbb{R}^d)$, and define the measurable function $\alpha : C([0, \infty[; \mathbb{R}^d) \times C([0, \infty[; \mathbb{R}^d) \to D(\mathbb{R}; \mathbb{R}^d)$ as

$$\alpha(\eta_-, \eta_+)(t) = \begin{cases} \eta_-(-t), & \text{if } t < 0, \\ \eta_+(t), & \text{if } t \geq 0. \end{cases}$$

Then $\tilde{\varphi} = (\mathbf{W} \otimes \mathbf{W}) \circ \alpha^{-1}$. Note that $\tilde{\varphi}$ is invariant under the shift $\tilde{\Theta}^s$.

Let $h : \tilde{E} \to \mathbb{R}^d$ be a continuous function (the observation function), so that $t \mapsto h(\xi_t)$ is càdlàg. By [22], we may define an $\overset{\frown}{\tilde{\mathcal{F}}}_{[s,t]}$-measurable map $\tilde{\Sigma}_{s,t}$ so that

$$\tilde{\Sigma}_{s,t}(\xi, \eta) = \exp\left( \int_s^t h(\xi_r) \cdot d\eta_r - \frac{1}{2} \int_s^t \|h(\xi_r)\|^2 \, dr \right) \qquad \text{for } \tilde{\varphi}\text{-a.e. } \eta \in \tilde{\Omega}^\eta$$

for every $\xi \in \tilde{\Omega}^\xi$. Note that $\tilde{\Sigma}_{s,t}$ is strictly positive by construction. We now define for every $s \leq t$ the probability kernel $\tilde{\Phi}_{s,t} : \tilde{\Omega}^\xi \times \overset{\frown}{\tilde{\mathcal{F}}}{}^\eta_{[s,t]} \to [0,1]$ as

$$\tilde{\Phi}_{s,t}(\xi, A) = \int I_A(\eta) \tilde{\Sigma}_{s,t}(\xi, \eta) \tilde{\varphi}(d\eta) \qquad \text{for all } A \in \overset{\frown}{\tilde{\mathcal{F}}}{}^\eta_{[s,t]}, \xi \in \tilde{\Omega}^\xi.$$

Define the process

$$\bar{\eta}_r = \eta_{r+s} - \eta_s - \int_s^{r+s} h(\xi_u) \, du.$$

Then by Girsanov's theorem, $(\bar{\eta}_r)_{r \in [0, t-s]}$ is a standard $d$-dimensional Wiener process under $\tilde{\Phi}_{s,t}(\xi, \cdot)$ for every $\xi \in \tilde{\Omega}^\xi$, as $t \mapsto h(\xi_t)$ is càdlàg and hence locally bounded (the usual conditions, which we have not assumed, are not needed for this to hold; see [38], Chapter 10). It remains to note that $\{\tilde{\Phi}_{s,t}(\xi, \cdot) : s \leq t\}$ is a consistent family, so there exists a probability kernel $\tilde{\Phi} : \tilde{\Omega}^\xi \times \overset{\frown}{\tilde{\mathcal{F}}}{}^\eta \to [0,1]$ with

$$\tilde{\Phi}(\xi, A) = \tilde{\Phi}_{s,t}(\xi, A) \qquad \text{for all } A \in \overset{\frown}{\tilde{\mathcal{F}}}{}^\eta_{[s,t]}, \xi \in \tilde{\Omega}^\xi, s \leq t,$$

by the usual Kolmogorov extension argument. It is easily verified that $\tilde{\Phi}$ satisfies the required properties, and Assumption 6.2 holds true by construction.

From this point onward we consider again the general continuous time setting (i.e., we do not assume white noise observations). The goal of this section is to extend several of our discrete time results to the continuous time setting. To this end, we will first prove the following counterpart of Theorem 4.2.

THEOREM 6.4.    *Suppose that Assumptions 6.1 and 6.2 are in force. Then*

$$\bigcap_{t \geq 0} \tilde{\mathscr{F}}_+^\eta \vee \sigma\{\xi_s : s \geq t\} = \tilde{\mathscr{F}}_+^\eta \quad \text{and} \quad \bigcap_{t \geq 0} \tilde{\mathscr{F}}_0^\eta \vee \tilde{\mathscr{F}}_{-t}^\xi = \tilde{\mathscr{F}}_0^\eta \qquad \mathbf{P}\text{-}a.s.$$

We now turn to the filter stability problem. As in discrete time, we must choose suitable versions of the regular conditional probabilities $\tilde{\mathbf{P}}^\mu(\xi_t \in \cdot | \tilde{\mathscr{F}}_{[0,t]}^\eta)$.

LEMMA 6.5.    *Suppose Assumption 6.2 holds. For any probability measure $\mu$ on $\mathcal{B}(\tilde{E})$, define a family of probability kernels $\tilde{\Pi}_t^\mu : \tilde{\Omega}^\eta \times \mathcal{B}(\tilde{E}) \to [0,1]$ $(t \geq 0)$ by*

$$\tilde{\Pi}_t^\mu(\eta, A) = \frac{\int I_A(\xi(t)) \tilde{\Sigma}_{0,t}(\xi, \eta) \tilde{\mathbf{P}}^\mu(d\xi)}{\int \tilde{\Sigma}_{0,t}(\xi, \eta) \tilde{\mathbf{P}}^\mu(d\xi)}.$$

*Then $\tilde{\Pi}_t^\mu$ is a version of the regular conditional probability $\tilde{\mathbf{P}}^\mu(\xi_t \in \cdot | \tilde{\mathscr{F}}_{[0,t]}^\eta)$.*

PROOF.    Apply the Bayes formula as in Lemma 5.1.    □

We can now prove a counterpart of Theorem 5.2. Note that the continuous time result yields a slightly weaker type of convergence than its discrete time counterpart; the reason for this choice is explained in the remark below.

THEOREM 6.6.    *Suppose that Assumptions 6.1 and 6.2 hold. Then*

$$\tilde{\mathbf{E}}^\mu(\|\tilde{\Pi}_t^\mu - \tilde{\Pi}_t^{\tilde{\pi}}\|_{\mathrm{TV}}) \xrightarrow{t \to \infty} 0 \quad \text{iff} \quad \|\tilde{\mathbf{P}}^\mu(\xi_t \in \cdot) - \tilde{\pi}\|_{\mathrm{TV}} \xrightarrow{t \to \infty} 0.$$

*Moreover, if*

$$\|\tilde{\mathbf{P}}^\mu(\xi_t \in \cdot) - \tilde{\pi}\|_{\mathrm{TV}} \xrightarrow{t \to \infty} 0 \quad \text{and} \quad \|\tilde{\mathbf{P}}^\nu(\xi_t \in \cdot) - \tilde{\pi}\|_{\mathrm{TV}} \xrightarrow{t \to \infty} 0,$$

*then $\tilde{\mathbf{E}}^\nu(\|\tilde{\Pi}_t^\mu - \tilde{\Pi}_t^{\tilde{\pi}}\|_{\mathrm{TV}}) \to 0$ as $t \to \infty$.*

REMARK 6.7.    Theorem 5.2 yields almost sure convergence of the filtering error, while Theorem 6.6 only gives convergence in $L^1$. The subtlety lies in the fact that convergence results for stochastic processes in continuous time, such as the martingale convergence theorem, require the choice of a modification of the stochastic process with appropriate continuity properties, and this typically requires that the filtrations satisfy the usual conditions (the associated $\sigma$-fields are therefore no longer countably generated). Though it seems very likely that such issues can be resolved with sufficient care, for example, along the lines of [39], we have chosen the simpler route which avoids unnecessary complications at the expense of a slightly weaker notion of convergence.

The remainder of this section is devoted to the proofs of Theorems 6.4 and 6.6.

6.2. *Reduction to discrete time.* The proofs in the continuous time setting can largely be reduced to our previous discrete time results. To this end, we begin by constructing a discrete time hidden Markov model, as defined in Section 3.1, which coincides with the continuous time model of this section.

The signal and observation state spaces for our discrete model are taken to be $E = D([0,1]; \tilde{E})$ and $F = D([0,1]; \tilde{F})$, respectively (recall that these Skorokhod spaces are themselves Polish). For the discrete time signal we will choose the $E$-valued process $X_n = (\xi_t)_{n \leq t \leq n+1}$, while we choose for the discrete time observations the $F$-valued process $Y_n = (\eta_t - \eta_n)_{n \leq t \leq n+1}$. We claim that these processes define a hidden Markov model in the sense of Section 3.1. Indeed, it is easily seen that $X_n$ is a Markov process with transition probability kernel

$$P(\xi', A) = \tilde{\mathbf{Q}}^{\xi'(1)}((\xi_t)_{0 \leq t \leq 1} \in A) \qquad \text{for all } \xi' \in E, A \in \mathcal{B}(E)$$

and invariant measure

$$\pi(A) = \tilde{\mathbf{P}}((\xi_t)_{0 \leq t \leq 1} \in A) \qquad \text{for all } A \in \mathcal{B}(E).$$

On the other hand, given $\tilde{\overset{\leftrightarrow}{\mathcal{F}}}{}^{\xi} = \mathcal{F}^X$, the random variables $Y_n$ are independent (as $\eta_t$ has conditionally independent increments given $\tilde{\overset{\leftrightarrow}{\mathcal{F}}}{}^{\xi}$) and we may define

$$\Phi((\xi(t))_{0 \leq t \leq 1}, A) = \tilde{\Phi}(\xi, Y_0 \in A) \qquad \text{for all } \xi \in \tilde{\Omega}^{\xi}, A \in \mathcal{B}(F),$$

where we have used that $\tilde{\Phi}(\xi, A)$ is $\tilde{\overset{\leftrightarrow}{\mathcal{F}}}{}^{\xi}_{[0,1]}$-measurable for $A \in \tilde{\overset{\leftrightarrow}{\mathcal{F}}}{}^{\eta}_{[0,1]}$ and that

$$\tilde{\mathbf{P}}(Y_n \in A | \tilde{\overset{\leftrightarrow}{\mathcal{F}}}{}^{\xi}) = \tilde{\Phi}(\xi, Y_n \in A) = \tilde{\Phi}(\tilde{\Theta}^n \xi, Y_0 \in A) = \Phi(X_n, A).$$

Having defined the kernels $P$ and $\Phi$ and the measure $\pi$, we may now construct the process $(X_n, Y_n)_{n \in \mathbb{Z}}$ on its canonical path space as in Section 3.1, and it is easily verified that the measures $\mathbf{P}$ and $\mathbf{P}^{\tilde{\mu}}$ coincide with the law of the process $(X_n, Y_n)$ under $\tilde{\mathbf{P}}$ and $\tilde{\mathbf{P}}^{\mu}$, respectively, where $\tilde{\mu} = \tilde{\mathbf{P}}^{\mu}(X_0 \in \cdot)$.

LEMMA 6.8. *Assumption 6.1 implies Assumption 3.1 for the discrete chain. Similarly, Assumption 6.2 implies Assumption 3.2 for the discrete chain.*

PROOF. By the Markov property, we find that

$$\|\tilde{\mathbf{Q}}^z((\xi_t)_{n \leq t \leq n+1} \in \cdot) - \pi\|_{\mathrm{TV}} = \|\tilde{\mathbf{Q}}^z(\xi_n \in \cdot) - \tilde{\pi}\|_{\mathrm{TV}}.$$

But note also that

$$\tilde{\mathbf{Q}}^z((\xi_t)_{n \leq t \leq n+1} \in \cdot) = \mathbf{P}^{\xi'}(X_{n+1} \in \cdot) \qquad \text{for all } \xi' \in E \text{ with } \xi'(1) = z.$$

The first statement follows directly. To prove the second statement, it suffices to note that under Assumption 6.2 we can write for $\xi \in \tilde{\Omega}^{\xi}$

$$\Phi((\xi_t)_{0 \leq t \leq 1}, A) = \int I_A((\eta_t - \eta_0)_{0 \leq t \leq 1}) \tilde{\Sigma}_{0,1}((\xi_t)_{0 \leq t \leq 1}, (\eta_t - \eta_0)_{0 \leq t \leq 1}) \tilde{\varphi}(d\eta),$$

so we may set $\varphi(A) = \tilde{\varphi}(Y_0 \in A)$ and $g(z, u) = \tilde{\Sigma}_{0,1}(z, u)$.  $\square$

The proof of Theorem 6.4 now follows immediately.

PROOF OF THEOREM 6.4. The result follows immediately from Theorem 4.2 in view of the fact that the measures $\tilde{\mathbf{P}}$ and $\mathbf{P}$ coincide.  $\square$

Before we proceed, let us prove a continuous time counterpart of Lemma 3.7.

LEMMA 6.9. *Suppose Assumption 6.2 holds. Let $\nu, \bar{\nu}$ be probability measures such that $\|\tilde{\mathbf{P}}^\nu(\xi_t \in \cdot) - \tilde{\mathbf{P}}^{\bar{\nu}}(\xi_t \in \cdot)\|_{\mathrm{TV}} \xrightarrow{t \to \infty} 0$. Then $\tilde{\mathbf{P}}^\nu|_{\widetilde{\mathscr{P}}_+^\eta} \sim \tilde{\mathbf{P}}^{\bar{\nu}}|_{\widetilde{\mathscr{P}}_+^\eta}$.*

PROOF. The result follows from Lemma 3.7, in view of the equivalence of the measures $\tilde{\mathbf{P}}^\mu$ and $\mathbf{P}^{\tilde{\mu}}$ ($\tilde{\mu} = \tilde{\mathbf{P}}^\mu(X_0 \in \cdot)$) for any $\mu$, using the same argument as in the proof of the first assertion of Lemma 6.8.  $\square$

6.3. *Proof of Theorem 6.6.* As in the discrete time setting, we begin by proving the only if part of Theorem 6.6. The proof is essentially identical.

LEMMA 6.10. *Suppose that Assumptions 6.1 and 6.2 hold and that*

$$\limsup_{t \to \infty} \|\tilde{\mathbf{P}}^\mu(\xi_t \in \cdot) - \tilde{\pi}\|_{\mathrm{TV}} > 0.$$

*Then we must have*

$$\limsup_{t \to \infty} \tilde{\mathbf{E}}^\mu(\|\tilde{\Pi}_t^\mu - \tilde{\Pi}_t^{\tilde{\pi}}\|_{\mathrm{TV}}) > 0.$$

PROOF. Let $\tilde{\mathbf{P}}^\mu(\xi_n \in \cdot) = \mu_n + \mu_n^\perp$ be the Lebesgue decomposition of $\tilde{\mathbf{P}}^\mu(\xi_n \in \cdot)$ with respect to $\tilde{\pi}$. In particular, $\mu_n \ll \tilde{\pi}$ and $\mu_n^\perp \perp \tilde{\pi}$, and there exists a set $S_n$ such that $\tilde{\pi}(S_n) = 0$ and $\mu_n^\perp(S_n^c) = 0$. We claim that

$$\limsup_{t \to \infty} \|\tilde{\mathbf{P}}^\mu(\xi_t \in \cdot) - \tilde{\pi}\|_{\mathrm{TV}} > 0 \quad \Longrightarrow \quad \limsup_{n \to \infty} \tilde{\mathbf{P}}^\mu(\xi_n \in S_n) > 0.$$

To see this, note that $(\xi_n)_{n \in \mathbb{Z}_+}$ is a discrete time Markov chain on the state space $\tilde{E}$. By [28], Theorem 7.2, Assumption 6.1 and $\tilde{\mathbf{P}}^\mu(\xi_n \in S_n) \to 0$ as $n \to \infty$ would imply that $\|\tilde{\mathbf{P}}^\mu(\xi_n \in \cdot) - \tilde{\pi}\|_{\mathrm{TV}} \to 0$ as $n \to \infty$. But $\|\tilde{\mathbf{P}}^\mu(\xi_t \in \cdot) - \tilde{\pi}\|_{\mathrm{TV}}$ is nonincreasing with $t$, so the latter implies that $\|\tilde{\mathbf{P}}^\mu(\xi_t \in \cdot) - \tilde{\pi}\|_{\mathrm{TV}} \to 0$ as $t \to \infty$. The claim is therefore established by contradiction.

Now note that it is easily established, using the expression in the proof of Lemma 6.5, that Assumption 6.2 implies $\tilde{\Pi}_{\underline{g}}^{\tilde{\pi}}(\eta, \cdot) \sim \tilde{\pi}$ for every $\eta \in \tilde{\Omega}^\eta$. Therefore, evidently $\tilde{\Pi}_n^{\tilde{\pi}}(\eta, S_n) = 0$ for all $\eta \in \tilde{\Omega}^\eta$, and we can estimate as follows:

$$\tilde{\Pi}_n^\mu(\eta, S_n) = |\tilde{\Pi}_n^\mu(\eta, S_n) - \tilde{\Pi}_n^{\tilde{\pi}}(\eta, S_n)| \le \|\tilde{\Pi}_n^\mu(\eta, \cdot) - \tilde{\Pi}_n^{\tilde{\pi}}(\eta, \cdot)\|_{\mathrm{TV}}.$$

In particular, we find that

$$\tilde{\mathbf{P}}^{\mu}(X_n \in S_n) = \tilde{\mathbf{E}}^{\mu}(\tilde{\Pi}_n^{\mu}((\eta_t)_{0 \leq t \leq n}, S_n)) \leq \tilde{\mathbf{E}}^{\mu}(\|\tilde{\Pi}_n^{\mu} - \tilde{\Pi}_n^{\tilde{\pi}}\|_{\mathrm{TV}})$$

and we must therefore have

$$\limsup_{n \to \infty} \tilde{\mathbf{E}}^{\mu}(\|\tilde{\Pi}_n^{\mu} - \tilde{\Pi}_n^{\tilde{\pi}}\|_{\mathrm{TV}}) > 0.$$

The proof is easily completed. $\square$

We now proceed to prove the converse assertion. One could attempt to adapt the corresponding discrete time proof to the current setting, but here we choose a different approach. First, we will show using Theorem 5.2 that

$$\|\tilde{\mathbf{P}}^{\mu}(\xi_t \in \cdot) - \tilde{\pi}\|_{\mathrm{TV}} \xrightarrow{t \to \infty} 0 \quad \text{and} \quad \|\tilde{\mathbf{P}}^{\nu}(\xi_t \in \cdot) - \tilde{\pi}\|_{\mathrm{TV}} \xrightarrow{t \to \infty} 0$$

implies that

$$\tilde{\mathbf{E}}^{\nu}(\|\tilde{\Pi}_n^{\mu} - \tilde{\Pi}_n^{\tilde{\pi}}\|_{\mathrm{TV}}) \xrightarrow{n \to \infty} 0,$$

where the limit as $n \to \infty$ is taken along the integers $n \in \mathbb{N}$. In the second step, we will show that the function

$$t \mapsto \tilde{\mathbf{E}}^{\nu}(\|\tilde{\Pi}_t^{\mu} - \tilde{\Pi}_t^{\tilde{\pi}}\|_{\mathrm{TV}})(t \in \mathbb{R}_+)$$

converges to a limit when we let $t \to \infty$ along the positive reals. Taken together, these two facts complete the proof of Theorem 6.6.

LEMMA 6.11. *Suppose that Assumptions 6.1 and 6.2 hold and that*

$$\|\tilde{\mathbf{P}}^{\mu}(\xi_t \in \cdot) - \tilde{\pi}\|_{\mathrm{TV}} \xrightarrow{t \to \infty} 0 \quad \text{and} \quad \|\tilde{\mathbf{P}}^{\nu}(\xi_t \in \cdot) - \tilde{\pi}\|_{\mathrm{TV}} \xrightarrow{t \to \infty} 0.$$

*Then* $\tilde{\mathbf{E}}^{\nu}(\|\tilde{\Pi}_n^{\mu} - \tilde{\Pi}_n^{\tilde{\pi}}\|_{\mathrm{TV}}) \xrightarrow{n \to \infty} 0 \ (n \in \mathbb{N})$.

PROOF. Let $\Pi_n^{\tilde{\mu}}$ and $\Pi_n^{\pi}$ be the filters for the discrete time chain as defined in Lemma 5.1, where $\tilde{\mu} = \tilde{\mathbf{P}}^{\mu}(X_0 \in \cdot)$. Note that, using the Markov property, we find that the condition of the current result implies that

$$\|\mathbf{P}^{\tilde{\mu}}(X_n \in \cdot) - \pi\|_{\mathrm{TV}} \xrightarrow{n \to \infty} 0.$$

Therefore, by Assumptions 6.1 and 6.2, Lemma 6.8 and Theorem 5.2, we find that

$$\|\Pi_n^{\tilde{\mu}} - \Pi_n^{\pi}\|_{\mathrm{TV}} \xrightarrow{t \to \infty} 0 \qquad \mathbf{P}^{\tilde{\mu}}\text{-a.s.}$$

It follows directly that

$$\|\Pi_n^{\tilde{\mu}}(Y, \xi(1) \in \cdot) - \Pi_n^{\pi}(Y, \xi(1) \in \cdot)\|_{\mathrm{TV}} \xrightarrow{n \to \infty} 0 \qquad \mathbf{P}^{\tilde{\mu}}\text{-a.s.}$$

But note that $\Pi_n^{\tilde{\mu}}(y, \xi(1) \in \cdot)$ and $\Pi_n^{\pi}(y, \xi(1) \in \cdot)$ are versions of the regular conditional probabilities

$$\tilde{\mathbf{P}}^{\mu}(\xi_{n+1} \in \cdot \,|\, \tilde{\mathscr{F}}_{[0,n+1]}^{\eta}) \quad \text{and} \quad \tilde{\mathbf{P}}(\xi_{n+1} \in \cdot \,|\, \tilde{\mathscr{F}}_{[0,n+1]}^{\eta}),$$

respectively. By the a.s. uniqueness of regular conditional probabilities and using Lemma 6.9 (which holds by virtue of Assumption 6.2), we therefore find that

$$\|\tilde{\Pi}_n^{\mu} - \tilde{\Pi}_n^{\pi}\|_{\mathrm{TV}} \xrightarrow{n \to \infty} 0 \qquad \tilde{\mathbf{P}}^{\nu}\text{-a.s.}$$

The result follows by dominated convergence. $\quad\square$

LEMMA 6.12. *Suppose that Assumption 6.2 holds and that*

$$\|\tilde{\mathbf{P}}^{\mu}(\xi_t \in \cdot) - \tilde{\pi}\|_{\mathrm{TV}} \xrightarrow{t \to \infty} 0 \quad \text{and} \quad \|\tilde{\mathbf{P}}^{\nu}(\xi_t \in \cdot) - \tilde{\pi}\|_{\mathrm{TV}} \xrightarrow{t \to \infty} 0.$$

*Then* $\tilde{\mathbf{E}}^{\nu}(\|\tilde{\Pi}_t^{\mu} - \tilde{\Pi}_t^{\pi}\|_{\mathrm{TV}})$ *is convergent as* $t \to \infty$ *($t \in \mathbb{R}_+$).*

PROOF. Let $\rho = (\mu + \tilde{\pi})/2$. Then we can establish, exactly as in the proof of Lemma 5.6, that we have $\tilde{\Pi}_t^{\mu} \ll \tilde{\Pi}_t^{\rho}$ and $\tilde{\Pi}_t^{\tilde{\pi}} \ll \tilde{\Pi}_t^{\rho}$ with

$$\frac{d\tilde{\Pi}_t^{\mu}}{d\tilde{\Pi}_t^{\rho}} = \frac{\tilde{\mathbf{E}}^{\rho}((d\mu/d\rho)(\xi_0)|\tilde{\mathscr{F}}_+^{\eta} \vee \tilde{\mathscr{F}}_{[t,\infty[}^{\xi})}{\tilde{\mathbf{E}}^{\rho}((d\mu/d\rho)(\xi_0)|\tilde{\mathscr{F}}_{[0,t]}^{\eta})},$$

$$\frac{d\tilde{\Pi}_t^{\tilde{\pi}}}{d\tilde{\Pi}_t^{\rho}} = \frac{\tilde{\mathbf{E}}^{\rho}((d\tilde{\pi}/d\rho)(\xi_0)|\tilde{\mathscr{F}}_+^{\eta} \vee \tilde{\mathscr{F}}_{[t,\infty[}^{\xi})}{\tilde{\mathbf{E}}^{\rho}((d\tilde{\pi}/d\rho)(\xi_0)|\tilde{\mathscr{F}}_{[0,t]}^{\eta})}, \qquad \tilde{\mathbf{P}}^{\rho}\text{-a.s.}$$

Note that $\tilde{\mathbf{E}}^{\rho}(d\tilde{\Pi}_t^{\mu}/d\tilde{\Pi}_t^{\rho}) = \tilde{\mathbf{E}}^{\rho}(d\tilde{\Pi}_t^{\tilde{\pi}}/d\tilde{\Pi}_t^{\rho}) = 1$ for all $t$. Now fix an arbitrary sequence $t_k \nearrow \infty$. By the martingale convergence theorem, we have $\tilde{\mathbf{P}}^{\rho}$-a.s.

$$\tilde{\mathbf{E}}^{\rho}\left(\frac{d\mu}{d\rho}(\xi_0)\,\Big|\,\tilde{\mathscr{F}}_{[0,t_k]}^{\eta}\right) \to \tilde{\mathbf{E}}^{\rho}\left(\frac{d\mu}{d\rho}(\xi_0)\,\Big|\,\tilde{\mathscr{F}}_+^{\eta}\right),$$

$$\tilde{\mathbf{E}}^{\rho}\left(\frac{d\tilde{\pi}}{d\rho}(\xi_0)\,\Big|\,\tilde{\mathscr{F}}_{[0,t_k]}^{\eta}\right) \to \tilde{\mathbf{E}}^{\rho}\left(\frac{d\tilde{\pi}}{d\rho}(\xi_0)\,\Big|\,\tilde{\mathscr{F}}_+^{\eta}\right)$$

as $k \to \infty$. Moreover, these quantities are $\tilde{\mathbf{P}}^{\rho}$-a.s. strictly positive by Lemma 6.9. Applying again the martingale convergence theorem, we find that $M_k^{\mu} := d\tilde{\Pi}_{t_k}^{\mu}/d\tilde{\Pi}_{t_k}^{\rho}$ and $M_k^{\tilde{\pi}} := d\tilde{\Pi}_{t_k}^{\tilde{\pi}}/d\tilde{\Pi}_{t_k}^{\rho}$ converge $\tilde{\mathbf{P}}^{\rho}$-a.s. to the random variables

$$M^{\mu} = \frac{\tilde{\mathbf{E}}^{\rho}((d\mu/d\rho)(\xi_0)|\bigcap_t \tilde{\mathscr{F}}_+^{\eta} \vee \tilde{\mathscr{F}}_{[t,\infty[}^{\xi})}{\tilde{\mathbf{E}}^{\rho}((d\mu/d\rho)(\xi_0)|\tilde{\mathscr{F}}_+^{\eta})}$$

and

$$M^{\tilde{\pi}} = \frac{\tilde{\mathbf{E}}^{\rho}((d\tilde{\pi}/d\rho)(\xi_0)|\bigcap_t \tilde{\mathscr{F}}_+^{\eta} \vee \tilde{\mathscr{F}}_{[t,\infty[}^{\xi})}{\tilde{\mathbf{E}}^{\rho}((d\tilde{\pi}/d\rho)(\xi_0)|\tilde{\mathscr{F}}_+^{\eta})},$$

respectively. Moreover, by the tower property of the conditional expectation, we have $\tilde{\mathbf{E}}^\rho(M^\mu|\widetilde{\mathscr{F}}_+^\eta) = 1$ and $\tilde{\mathbf{E}}^\rho(M^{\tilde{\pi}}|\widetilde{\mathscr{F}}_+^\eta) = 1$ $\tilde{\mathbf{P}}^\rho$-a.s. Therefore, $\tilde{\mathbf{E}}^\rho(M^\mu) = \tilde{\mathbf{E}}^\rho(M^{\tilde{\pi}}) = 1$, so that $M_k^\mu \to M^\mu$ and $M_k^{\tilde{\pi}} \to M^{\tilde{\pi}}$ in $L^1(\tilde{\mathbf{P}}^\rho)$ by Scheffé's lemma.

Let us write, for simplicity, $N_k = |M_k^\mu - M_k^{\tilde{\pi}}|$ and $N = |M^\mu - M^{\tilde{\pi}}|$. Then

$$\tilde{\mathbf{E}}^\rho(|\tilde{\mathbf{E}}^\rho(N_k|\widetilde{\mathscr{F}}_{[0,t_k]}^\eta) - \tilde{\mathbf{E}}^\rho(N|\widetilde{\mathscr{F}}_+^\eta)|)$$

$$\leq \tilde{\mathbf{E}}^\rho(|\tilde{\mathbf{E}}^\rho(N_k|\widetilde{\mathscr{F}}_{[0,t_k]}^\eta) - \tilde{\mathbf{E}}^\rho(N|\widetilde{\mathscr{F}}_{[0,t_k]}^\eta)|) + \tilde{\mathbf{E}}^\rho(|\tilde{\mathbf{E}}^\rho(N|\widetilde{\mathscr{F}}_{[0,t_k]}^\eta) - \tilde{\mathbf{E}}^\rho(N|\widetilde{\mathscr{F}}_+^\eta)|)$$

$$\leq \tilde{\mathbf{E}}^\rho(|N_k - N|) + \tilde{\mathbf{E}}^\rho(|\tilde{\mathbf{E}}^\rho(N|\widetilde{\mathscr{F}}_{[0,t_k]}^\eta) - \tilde{\mathbf{E}}^\rho(N|\widetilde{\mathscr{F}}_+^\eta)|)$$

$$\leq \tilde{\mathbf{E}}^\rho(|M_k^\mu - M_k^{\tilde{\pi}} - M^\mu + M^{\tilde{\pi}}|) + \tilde{\mathbf{E}}^\rho(|\tilde{\mathbf{E}}^\rho(N|\widetilde{\mathscr{F}}_{[0,t_k]}^\eta) - \tilde{\mathbf{E}}^\rho(N|\widetilde{\mathscr{F}}_+^\eta)|)$$

$$\leq \tilde{\mathbf{E}}^\rho(|M_k^\mu - M^\mu|) + \tilde{\mathbf{E}}^\rho(|M_k^{\tilde{\pi}} - M^{\tilde{\pi}}|) + \tilde{\mathbf{E}}^\rho(|\tilde{\mathbf{E}}^\rho(N|\widetilde{\mathscr{F}}_{[0,t_k]}^\eta) - \tilde{\mathbf{E}}^\rho(N|\widetilde{\mathscr{F}}_+^\eta)|),$$

where we have used the inverse triangle inequality to establish that $|N_k - N| \leq |M_k^\mu - M_k^{\tilde{\pi}} - M^\mu + M^{\tilde{\pi}}|$. By the martingale convergence theorem and the convergence of $M_k^\mu$ and $M_k^{\tilde{\pi}}$, the right-hand side of this expression converges to zero as $k \to \infty$. But note that $\|\tilde{\Pi}_{t_k}^\mu - \tilde{\Pi}_{t_k}^{\tilde{\pi}}\|_{\mathrm{TV}} = \tilde{\mathbf{E}}^\rho(N_k|\widetilde{\mathscr{F}}_{[0,t_k]}^\eta)$ $\tilde{\mathbf{P}}^\rho$-a.s., so we have

$$\|\tilde{\Pi}_{t_k}^\mu - \tilde{\Pi}_{t_k}^{\tilde{\pi}}\|_{\mathrm{TV}} \xrightarrow{k \to \infty} \tilde{\mathbf{E}}^\rho(N|\widetilde{\mathscr{F}}_+^\eta) \qquad \text{in } L^1(\tilde{\mathbf{P}}^\rho).$$

In particular, $\|\tilde{\Pi}_{t_k}^\mu - \tilde{\Pi}_{t_k}^{\tilde{\pi}}\|_{\mathrm{TV}}$ converges to $\tilde{\mathbf{E}}^\rho(N|\widetilde{\mathscr{F}}_+^\eta)$ in $\tilde{\mathbf{P}}^\rho$-probability. But

$$\|\tilde{\mathbf{P}}^\nu(\xi_t \in \cdot) - \tilde{\mathbf{P}}^\rho(\xi_t \in \cdot)\|_{\mathrm{TV}}$$

$$\leq \tfrac{1}{2}(\|\tilde{\mathbf{P}}^\nu(\xi_t \in \cdot) - \tilde{\mathbf{P}}^\mu(\xi_t \in \cdot)\|_{\mathrm{TV}} + \|\tilde{\mathbf{P}}^\nu(\xi_t \in \cdot) - \tilde{\pi}\|_{\mathrm{TV}})$$

$$\leq \tfrac{1}{2}(\|\tilde{\mathbf{P}}^\mu(\xi_t \in \cdot) - \tilde{\pi}\|_{\mathrm{TV}} + 2\|\tilde{\mathbf{P}}^\nu(\xi_t \in \cdot) - \tilde{\pi}\|_{\mathrm{TV}}) \xrightarrow{t \to \infty} 0,$$

so by Lemma 6.9 we find that $\|\tilde{\Pi}_{t_k}^\mu - \tilde{\Pi}_{t_k}^{\tilde{\pi}}\|_{\mathrm{TV}}$ converges to $\tilde{\mathbf{E}}^\rho(N|\widetilde{\mathscr{F}}_+^\eta)$ in $\tilde{\mathbf{P}}^\nu$-probability. Thus, we have $\tilde{\mathbf{E}}^\nu(\|\tilde{\Pi}_{t_k}^\mu - \tilde{\Pi}_{t_k}^{\tilde{\pi}}\|_{\mathrm{TV}}) \to \tilde{\mathbf{E}}^\nu(\tilde{\mathbf{E}}^\rho(N|\widetilde{\mathscr{F}}_+^\eta))$ by dominated convergence. But as this holds for any sequence $t_k \nearrow \infty$, the result follows. $\square$

**7. On the result of Kunita and necessity of the ergodic condition.** In Sections 5 and 6 we explored the consequences of our main results for the stability of nonlinear filters. Our results also have implications for other asymptotic properties of the filter, however, in particular for the uniqueness of the invariant measure as studied in [23]. The aim of this section is to briefly outline the connection with [23], and to compare our assumptions to those made in the work of Kunita.

Kunita's original paper [23] investigated the continuous time setting with compact signal state space and white noise type observations. His approach has been extended to locally compact [32] and Polish [2] signal state spaces on the one hand, and to discrete time models on locally compact [32] and Polish [17] signal state spaces on the other hand. None of these papers resolve the gap in [23], however; we refer to [4] for further discussion and references. For simplicity and concreteness, we will restrict our discussion below to the original setting of Kunita. However, the results in this paper apply to all settings considered in the above references, and the reader interested in the ergodic properties of the nonlinear filter can directly read off the relevant results from these papers.

In [23], the signal process $(\xi_t)_{t\in\mathbb{R}}$ is a stationary, time-homogeneous Feller–Markov process on a compact Polish state space $\tilde{E}$ with stationary measure $\tilde{\pi}$ under $\tilde{\mathbf{P}}$, and the $\mathbb{R}^d$-valued observation process $(\eta_t)_{t\in\mathbb{R}}$ is defined as

$$\eta_t = \int_0^t h(\xi_s)\,ds + W_t,$$

where $(W_t)_{t\in\mathbb{R}}$ is a two-sided Wiener process and $h\colon \tilde{E}\to\mathbb{R}^d$ is a continuous function. Kunita establishes that the filter $\tilde{\Pi}_t^{\tilde{\pi}}$, when seen as a measure-valued random process, is itself a Feller–Markov process, and we are interested in the ergodic properties of this process. In particular, [23] yields the following statement.

PROPOSITION 7.1.  *There exists at least one invariant measure for the filter whose barycenter is $\tilde{\pi}$. If $\bigcap_{t\geq 0}\tilde{\mathscr{F}}_0^{\eta}\vee\tilde{\mathscr{F}}_{-t}^{\xi}=\tilde{\mathscr{F}}_0^{\eta}$ $\tilde{\mathbf{P}}$-a.s. holds true, then there is only one such invariant measure. If in addition $\tilde{\pi}$ is the unique invariant measure of the signal, then the invariant measure of the filter is unique.*

In [23], it is assumed that the $\tilde{\mathbf{P}}$-a.s. triviality of the tail $\sigma$-field $\bigcap_{t\geq 0}\tilde{\mathscr{F}}_{-t}^{\xi}$, or, equivalently [34], Proposition 3, the condition

$$\int |\tilde{\mathbf{P}}^{\delta_z}(\xi_t\in A) - \tilde{\pi}(A)|\,\tilde{\pi}(dz) \xrightarrow{t\to\infty} 0 \qquad \text{for all } A\in\mathcal{B}(\tilde{E}),$$

is already sufficient to establish $\bigcap_{t\geq 0}\tilde{\mathscr{F}}_0^{\eta}\vee\tilde{\mathscr{F}}_{-t}^{\xi}=\tilde{\mathscr{F}}_0^{\eta}$ $\tilde{\mathbf{P}}$-a.s. As we have argued before, however, this statement is not at all obvious. On the other hand, by the continuous time version of our main result (Theorem 6.4), it follows that

$$\int \sup_{A\in\mathcal{B}(\tilde{E})} |\tilde{\mathbf{P}}^{\delta_z}(\xi_t\in A) - \tilde{\pi}(A)|\,\tilde{\pi}(dz) \xrightarrow{t\to\infty} 0$$

does in fact guarantee that $\bigcap_{t\geq 0}\tilde{\mathscr{F}}_0^{\eta}\vee\tilde{\mathscr{F}}_{-t}^{\xi}=\tilde{\mathscr{F}}_0^{\eta}$ $\tilde{\mathbf{P}}$-a.s. [that this condition is equivalent to Assumption 6.1 follows from the fact that $\|\tilde{\mathbf{P}}^{\mu}(\xi_t\in\cdot) - \tilde{\pi}\|_{\mathrm{TV}}$ is nonincreasing]. This condition covers most, but not all, of the models that satisfy Kunita's condition, and we have thus partially resolved the gap in his proof. Whether Kunita's condition is already sufficient to guarantee uniqueness of the invariant measure with barycenter $\tilde{\pi}$ remains an open problem.

Besides sufficiency of the ergodic condition, it is interesting to ask whether such a condition is necessary for uniqueness of the invariant measure. Theorem 3.3 of [23] states that Kunita's condition is in fact necessary for uniqueness of the invariant measure with barycenter $\tilde{\pi}$, but this does not appear to be correct. As the following example shows, neither our condition nor Kunita's condition is necessary.

EXAMPLE 7.2. Consider the signal on $\tilde{E} = [0, 1]$ such that $\xi_t = \xi_0$ for all $t \in \mathbb{R}$ $\tilde{\mathbf{P}}$-a.s., and let $\tilde{\pi}$ be the Lebesgue measure on $[0, 1]$. We choose the observation function $h(x) = x$. This model fits entirely within the current setting.

Let us first show that the signal does not satisfy Kunita's condition (and hence it does not satisfy our assumptions, which are stronger than Kunita's). Note that

$$\overset{\circ\circ}{\mathcal{F}}{}^{\xi}_{-t} = \sigma\{\xi_s : s \leq -t\} = \sigma\{\xi_0\} \qquad \tilde{\mathbf{P}}\text{-a.s. for all } t \in \mathbb{R}.$$

Therefore, $\tilde{\mathbf{P}}$-a.s. $\bigcap_{t \geq 0} \overset{\circ\circ}{\mathcal{F}}{}^{\xi}_{-t} = \sigma\{\xi_0\}$, which is certainly not $\tilde{\mathbf{P}}$-a.s. trivial.

We claim that nonetheless $\bigcap_{t \geq 0} \overset{\circ\circ}{\mathcal{F}}{}^{\eta}_{0} \vee \overset{\circ\circ}{\mathcal{F}}{}^{\xi}_{-t} = \overset{\circ\circ}{\mathcal{F}}{}^{\eta}_{0}$ $\tilde{\mathbf{P}}$-a.s., so the invariant measure of the filter with barycenter $\tilde{\pi}$ is unique. Clearly it suffices to show that

$$\overset{\circ\circ}{\mathcal{F}}{}^{\xi}_{-t} = \sigma\{\xi_0\} \subset \overset{\circ\circ}{\mathcal{F}}{}^{\eta}_{0} \qquad \tilde{\mathbf{P}}\text{-a.s.}$$

for all $t \geq 0$. But note that $\eta_t = \xi_0 t + W_t$ for all $t \in \mathbb{R}$, so

$$\limsup_{t \to -\infty} \frac{\eta_t}{t} = \xi_0 \qquad \tilde{\mathbf{P}}\text{-a.s.}$$

The claim is therefore established.

The previous example highlights a possibility which is not considered in this paper. Returning to our canonical model, suppose that the tail $\sigma$-field $\mathcal{T}^X$ is not $\mathbf{P}$-a.s. trivial (so the signal is not ergodic), but that $\mathcal{T}^X \subset \mathcal{F}^Y_{[0, \infty[}$ $\mathbf{P}$-a.s. Then, if it could somehow be established that the exchange of intersection and supremum is permitted, we would still obtain the identity

$$\bigcap_{n \geq 0} \mathcal{F}^Y_{[0, \infty[} \vee \mathcal{F}^X_{[n, \infty[} \overset{?}{=} \mathcal{F}^Y_{[0, \infty[} \vee \bigcap_{n \geq 0} \mathcal{F}^X_{[n, \infty[} = \mathcal{F}^Y_{[0, \infty[} \qquad \mathbf{P}\text{-a.s.},$$

and therefore also the associated implications for the stability properties and for the uniqueness of the invariant measure of the filter. The condition $\mathcal{T}^X \subset$ $\mathcal{F}^Y_{[0,\infty[}$ is closely related to the notion of detectability which is shown in [36] to be necessary and sufficient for the stability of the filter (in a suitable sense) for models with a finite signal state space and nondegenerate observations. Whether such a necessary and sufficient condition can be obtained for more general models in the absence of an ergodicity assumption is an interesting topic for further investigation.

DEPARTMENT OF OPERATIONS RESEARCH
AND FINANCIAL ENGINEERING
PRINCETON UNIVERSITY
PRINCETON, NEW JERSEY 08544
USA
E-MAIL: rvan@princeton.edu